\documentclass[a4paper,10pt]{article}
\usepackage{etex}
\usepackage{shuffle}
\usepackage{hyperref}
\usepackage{amsmath,amsfonts,amssymb}
\usepackage{fancyhdr}
\usepackage{mathrsfs,amssymb}
\usepackage{amsmath}
\usepackage{amssymb}
\usepackage{empheq}
\usepackage{caption}
\usepackage{enumerate}
\usepackage{tikz}
\usepackage{float}
\usepackage[a4paper]{geometry}
\usepackage[latin1]{inputenc}
\usepackage[T1]{fontenc}
\usepackage{lmodern}
\usepackage{amsthm}
\usepackage{dsfont}
\usepackage[nottoc]{tocbibind}
\usepackage{color}
\usepackage{fancyhdr,blindtext}
\usepackage[arrow, matrix, curve]{xy}
\usepackage{array}
\usepackage{stmaryrd}
\usepackage{bm}
 
\usepackage[vcentermath]{youngtab}
\newcolumntype{x}[1]{>{\centering\arraybackslash\hspace{0pt}}m{#1}}

\fancypagestyle{plain}{%
\fancyhf{}}

\renewcommand{\sectionmark}[1]%
{\markboth{#1}{}}

\makeindex

\theoremstyle{definition}
\newtheorem{ex}{\bfseries \upshape Example}[section]
\newtheorem{dfn}[ex]{\bfseries \upshape Definition}
\newtheorem{nota}[ex]{\bfseries \upshape Notation}
\newtheorem{rem}[ex]{\bfseries \upshape Remark}
\newtheorem{conj}[ex]{\bfseries \upshape Conjecture}
\newtheorem{prop}[ex]{\bfseries \upshape Proposition}
\newtheorem{lem}[ex]{\bfseries \upshape Lemma}
\newtheorem{thm}[ex]{\bfseries \upshape Theorem}

\newtheorem{constr}[ex]{\bfseries \upshape Construction}
\newtheorem{thmx}{Theorem}
\newtheorem{qu}{Question}

\theoremstyle{plain}

\newtheorem{cor}[ex]{\bfseries \upshape Corollary}

\newenvironment{prf}{\begin{proof}[{\bf Proof}]}{\end{proof}}
\newcommand{\h}{\mathfrak H}
\newcommand{\N}{\ensuremath{\mathds{N}}}	
\newcommand{\Z}{\ensuremath{\mathds{Z}}}	
\newcommand{\Q}{\ensuremath{\mathds{Q}}}	
\newcommand{\R}{\ensuremath{\mathbb{R}}}	
\newcommand{\C}{\ensuremath{\mathds{C}}}
\newcommand{\Ha}{\ensuremath{\mathbb{H}}}

\newcommand{\Li}{\operatorname{Li}}
\newcommand{\Lit}{\widetilde{\operatorname{Li}}}
\newcommand{\dif}{ \operatorname{d}_q}

\newcommand{\grw}{ \operatorname{gr}^{\operatorname{W}}}

\newcommand{\filw}{ \operatorname{Fil}^{\operatorname{W}}}
\newcommand{\fild}{ \operatorname{Fil}^{\operatorname{D}}}
\newcommand{\fille}{ \operatorname{Fil}^{\operatorname{L}} }
\newcommand{\filwle}{ \operatorname{Fil}^{\operatorname{W},\operatorname{L}}}

\DeclareMathOperator{\Sl}{SL}

\DeclareMathOperator{\MD}{\mathcal{MD}}

\DeclareMathOperator{\bMD}{\mathcal{BD}}
\DeclareMathOperator{\bMDG}{\bMD_{gen}}
\DeclareMathOperator{\MDA}{q\mathcal{M}\mathcal{Z}}

\DeclareMathOperator{\MZ}{\mathcal{MZ}}
\DeclareMathOperator{\MZB}{\mathcal{MZB}}

\newcommand{\lemref}[1] {Lemma \ref{#1}}

\numberwithin{equation}{section}

\newcommand{\lar}[1]{\overset{\text{\tiny$\bm\leftarrow$}}{#1}}
\newcommand\zetarev{\lar{\zeta}}
\newcommand\zetashrev{\lar{\zeta^\sh}}
\newcommand\zetastrev{\lar{\zeta^\ast}}
\newcommand\grev{\lar{G}}
\newcommand\gshrev{\lar{G^\sh}}
\newcommand\gstrev{\lar{G^\ast}}
\newcommand\gstmrev{\lar{G^{\ast,M}}}

\DeclareFontEncoding{OT2}{}{}
\DeclareFontSubstitution{OT2}{cmr}{m}{n}
\DeclareFontFamily{OT2}{cmr}{\hyphenchar\font45 }
\DeclareFontShape{OT2}{cmr}{m}{n}{<5><6><7><8><9>gen*wncyr<10><10.95><12><14.4><17.28><20.74><24.88>wncyr10}{}
\DeclareFontShape{OT2}{cmr}{b}{n}{<5><6><7><8><9>gen*wncyb<10><10.95><12><14.4><17.28><20.74><24.88>wncyb10}{}
\DeclareMathAlphabet{\mathcyr}{OT2}{cmr}{m}{n}
\DeclareMathAlphabet{\mathcyb}{OT2}{cmr}{b}{n}
\SetMathAlphabet{\mathcyr}{bold}{OT2}{cmr}{b}{n}
\newcommand{\sh}{\shuffle}
\newcommand{\mtt}[3] {{\Large \begin{vsmallmatrix} 
  #1\\
 	#2 \\
	#3
\end{vsmallmatrix}}}

\DeclareRobustCommand{\mb}{\genfrac{[}{]}{0pt}{}}
\DeclareRobustCommand{\mt}{\genfrac{|}{|}{0pt}{}}

\begin{document}
\title{{ \bf The algebra of bi-brackets and \\regularised multiple Eisenstein series}}
\author{{\sc Henrik Bachmann}}

\date{\today}
\maketitle
\abstract{We study the algebra of certain $q$-series, called bi-brackets, whose coefficients are given by weighted sums over partitions. These series incorporate the theory of modular forms for the full modular group as well as the theory of multiple zeta values (MZV) due to their appearance in the Fourier expansion of regularised multiple Eisenstein series. Using the conjugation of partitions we obtain linear relations between bi-brackets, called the partition relations, which yield naturally two different ways of expressing the product of two bi-brackets similar to the stuffle and shuffle product of multiple zeta values. Bi-brackets are generalizations of the generating functions of multiple divisor sums, called brackets,  $[s_1,\dots,s_l]$ studied in \cite{BK}. We use the algebraic structure of bi-brackets to define further $q$-series $[s_1,\dots,s_l]^\sh$ and $[s_1,\dots,s_l]^\ast$  which satisfy the shuffle and stuffle product formulas of MZV by using results about quasi-shuffle algebras introduced by Hoffman. In \cite{BT} regularised multiple Eisenstein series $G^\sh$ were defined, by using an explicit connection to the coproduct on formal iterated integrals. These satisfy the shuffle product formula. Applying the same concept for the coproduct on quasi-shuffle algebras enables us to define multiple Eisenstein series $G^\ast$ satisfying the stuffle product. We show that both $G^\sh$ and $G^\ast$ are given by linear combinations of products of MZV and bi-brackets. Comparing these two regularized multiple Eisenstein series enables us to obtain finite double shuffle relations for multiple Eisenstein series in low weights which extend the relations proven in \cite{BT}.}
\tableofcontents
\newpage

\section{Introduction}
Multiple zeta values are natural generalizations of the Riemann zeta values that are defined for integers $s_1 > 1$ and $s_i \geq 1$ for $i>1$ by
\[ \zeta( s_1 , \dots , s_l ) :=  \sum_{n_1 >n_2 >\dots>n_l>0 } \frac{1}{n_1^{s_1} \dots n_l^{s_l} } \,. \]
Because of its occurrence in various fields of mathematics and physics these real numbers 
are of particular interest.
The $\Q$-vector space of all multiple zeta values of weight $k$ is then given by 
\[ \MZ_k :=\big <\,\zeta(s_1,\dots,s_l) \, \big| \, s_1 + \dots + s_l = k 
\textrm{ and } l>0 \big>_\Q.  \]
It is well known that  the product of two multiple zeta values can be written as a linear combination of multiple zeta values of the same weight by using the stuffle or shuffle relations. Thus they generate a $\Q$-algebra $\MZ$. There are several connections of these numbers to modular forms for the full modular group. Some of them are treated in \cite{gkz}, where connections of double zeta values and modular forms are described. One of them is given by double Eisenstein series $G_{s_1,s_2} \in \C[[q]]$ which are the length two version of classical Eisenstein series and which are given by a double sum over ordered lattice points. These functions have a Fourier expansion given by sums of products of MZV and certain $q$-series with the double zeta value $\zeta(s_1,s_2)$ as their constant term. In \cite{Bach} the author treated the multiple case and calculated the Fourier expansion of multiple Eisenstein series (MES) $G_{s_1,\dots,s_l} \in \C[[q]]$. The result of \cite{Bach} was that the Fourier expansion of MES is again a linear combination of MZV and $q$-series $[s_1,\dots,s_l]\in \Q[[q]]$, called brackets, with the corresponding MZV as the constant term. For example it is
\begin{align*}
G_{3,2,2}(\tau) =& \zeta(3,2,2) + \left( \frac{54}{5} \zeta(2,3) + \frac{51}{5} \zeta(3,2) \right) (2 \pi i)^2 [2] + \frac{16}{3} \zeta  (2,2) (2 \pi  i)^3 [3] \\
&+3 \zeta(3) (2\pi i)^4 [2,2] + 4 \zeta(2) (2 \pi i)^5 [3,2] + (2 \pi i)^7 [3,2,2]\,.
\end{align*}
It turned out that the $q$-series $[s_1,\dots,s_l]$, whose coefficients $a_n$ are given by weighted sums over partitions of $n$, are, independently to their appearance in the Fourier expansion of MES, very interesting objects and therefore they were studied on their own in \cite{BK}. There the authors studied the algebraic structure of the space $\MD$ spanned by these brackets and we will refine, generalize and use some of the results in this note. 

Due to convergence issues the MES are just defined for $s_1,\dots,s_l \geq 2$ and therefore there are a lot more MZV than MES. A natural question was therefore the following
\begin{qu}\emph{
 What is a "good" definition of a "regularised" multiple Eisenstein series, such that for each multiple zeta value  $\zeta(s_1,\dots,s_l)$ with $s_1>1$,$s_2,\dots,s_l \geq 1$ there is a multiple Eisenstein series   
\[ G^{reg}_{s_1,\dots,s_l} =  \zeta(s_1,\dots,s_l)+ \sum_{n>0} a_n q^n \in \C[[q]] \]
with this multiple zeta values as the constant term in its Fourier expansion and which equals the original multiple Eisenstein series in the case $s_1,\dots,s_l \geq 2$?}
\end{qu}
 By "good" we mean 
that these multiple Eisenstein series should have the same, or at least as much as possible, algebraic structure as multiple zeta values, i.e. they should fulfill the shuffle or/and the stuffle product. In \cite{BT} the authors addressed this question and they define (shuffle) regularised MES $G^\sh_{s_1,\dots,s_l}$, defined for all $s_1,\dots,s_l \in \N$, which coincide with the $G_{s_1, \dots,s_l}$ in the case $s_1,\dots,s_l \geq 2$ and which fulfill the shuffle product. In their construction the authors consider certain $q$-series similar to the brackets which also fulfill the shuffle product.

In this note we want to consider a more general class of $q$-series which we call bi-brackets. We will see that the $q$-series appearing in the construction in \cite{BT} are linear combination of bi-brackets. Furthermore we will address the above question with respect to the stuffle product and we will construct another (stuffle) regularised type of MES, denoted by $G^\ast_{s_1,\dots,s_l}$, satisfying the stuffle product formula. The bi-brackets will also appear there and we will be able to write $G^\sh$ and $G^\ast$ as sums of products of MZV and bi-brackets which then enables us to compare these two types of regularised MES.

Even when one is not interested in the question of extending the definition of MES we want to emphasize the reader that these $q$-series are interesting by their own rights, since they give a $q$-analogue of multiple zeta values with a nice algebraic structure. These q-analogues have two ways to write the product of two such series similar to the shuffle and the stuffle product for MZV. 
For $s_1,\dots,s_l \geq 1$, $r_1,\dots,r_l \geq 0$ these $q$-series, which we call bi-brackets, are given by
\begin{align*}
\mb{s_1, \dots , s_l}{r_1,\dots,r_l} :=\sum_{\substack{u_1 > \dots > u_l> 0 \\ v_1, \dots , v_l >0}} \frac{u_1^{r_1}}{r_1!} \dots \frac{u_l^{r_l}}{r_1!} \cdot \frac{v_1^{s_1-1} \dots v_l^{s_l-1}}{(s_1-1)!\dots(s_l-1)!}   \cdot q^{u_1 v_1 + \dots + u_l v_l} \in \Q[[q]]\,.
\end{align*}
In the first section we will interpret this sum as a weighted sum over partitions of a natural number $n$. The conjugation of partitions will give us linear relations between the bi-brackets which we therefore call the partition relation. We use this relation to prove a stuffle and shuffle analogue of the product of two bi-brackets and obtain for example 
\begin{align*}
\mb{2,3}{0,0}+\mb{3,2}{0,0}+\mb{5}{0}-\frac{1}{12} \mb{3}{0} &= \mb{2}{0} \cdot \mb{3}{0} = \mb{2,3}{0,0}+3 \mb{3,2}{0,0}+6 \mb{4,1}{0,0}-3 \mb{4}{0}+3 \mb{4}{1} \,.
\end{align*}
Compare this with the "real" stuffle and shuffle product of multiple zeta values
\[ \zeta(2,3) + \zeta(3,2) + \zeta(5) = \zeta(2) \cdot \zeta(3) = \zeta(2,3) + 3 \zeta(3,2) + 6 \zeta(4,1) \,. \] 
Using the algebraic structure of the space of bi-brackets we define a shuffle $[s_1,\dots,s_l]^\sh$ and stuffle $[s_1,\dots,s_l]^\ast$ version of the ordinary brackets as certain linear combination of bi-brackets. These objects fulfill the same shuffle and stuffle products as multiple zeta values. Both constructions use the theory of quasi-shuffle algebras developed by Hoffman in \cite{H1}. We end the introduction by summarizing the results of this paper on bi-brackets and regularised multiple Eisenstein series in the following two vaguely formulated theorems: 
\begin{thmx}
\begin{enumerate}[i)]
\item The space $\bMD$ spanned by all bi-brackets $ \mb{s_1, \dots , s_l}{r_1,\dots,r_l}$  forms a $\Q$-algebra with the space of (quasi-)modular forms and the space  $\MD$ of brackets as subalgebras. There are two ways to express the product as a linear combination of bi-brackets which yields a large family of linear relations. 
\item There are two subalgebras $\MD^\sh \subset \bMD$ and $\MD^\ast \subset \MD$ spanned by elements $[s_1,\dots,s_l]^\sh$ and $[s_1,\dots,s_l]^\ast$ which fulfill the shuffle and stuffle products, respectively, and which are in the length one case given by the bracket $[s_1]$.
\end{enumerate}
\end{thmx}
For example we have similarly to the relation between MZV above
\[ [2,3]^\ast + [3,2]^\ast + [5] = [2] \cdot [3] = [2,3]^\sh + 3 [3,2]^\sh + 6 [4,1]^\sh \,. \] 
Denote by $\MZB \subset \C[[q]]$ the space of all formal power series in $q$ which can be written as a linear combination of products of MZV, powers of $(-2\pi i)$ and bi-brackets.
\begin{thmx}
\begin{enumerate}[i)]
\item  The shuffle regularised multiple Eisenstein series $G_{s_1,\dots,s_l}^\sh \in \C[[q]]$ defined in \cite{BT} can be written as a linear combination of products of MZV, powers of $(-2\pi i)$ and shuffle brackets $[r_1,\dots,r_m]^\sh$, i.e. they are elements of the space $\MZB$.
\item For all $s_1,\dots,s_l \in \N$ and $M \in \N$ there are $q$-series $G_{s_1,\dots,s_l}^{\ast,M} \in \C[[q]]$ (see Definition \ref{def:gast}) which fulfill the stuffle product. If the limit $G^\ast_{s_1,\dots,s_l} :=  \lim_{M\to\infty} G_{s_1,\dots,s_l}^{\ast,M}$ exists it will be an element in $\MZB$ which still fulfills the stuffle product. In that case the $q$-series $G^\ast_{s_1,\dots,s_l}$ will be called stuffle regularised multiple Eisenstein series.

\item For $s_1,\dots,s_l \geq 2$ both regularised multiple Eisenstein series equal the classical multiple Eisenstein series, i.e. we have
\[ G_{s_1,\dots,s_l} = G^\sh_{s_1,\dots,s_l} = G^\ast_{s_1,\dots,s_l} \,.\] 
\end{enumerate}
\end{thmx}

{\bf Content of this paper}: 
In section $2$ we will introduce bi-brackets and their generating series. We will show that there are a natural linear relations between bi-brackets, called the partition relations. In section $3$ we prove that the algebra of bi-brackets has the structure of a quasi-shuffle algebra in the sense of \cite{H1}. The partition relation will yield another way of multiplying two bi-brackets which differs from the quasi-shuffle product and which therefore yields linear relations similar to the double shuffle relations of MZV. The connection to modular forms and the derivatives of bi-brackets will be subject of section $4$. We will see that relations between bi-brackets can be used to prove relations between modular forms and vice versa. Section $5$ will be devoted to the definition of the brackets $[s_1,\dots,s_l]^\sh$ and $[s_1,\dots,s_l]^\ast$. For this we will recall the algebraic setup of Hoffman in this section. Finally in section $4$ we will recall the results of \cite{BT} and the definition of the shuffle regularised MES $G^\sh$. After this we will define the stuffle regularised MES $G^{\ast,M}$ and $G^{\ast}$ by using a similar approach as in the definition of $G^\sh$. We end section $4$ by comparing these two regularised MES in low weights. 

\subsubsection*{Acknowledgment}
The author would like to thank Herbert Gangl, Ulf Kühn, Koji Tasaka, Jianqiang Zhao and Wadim Zudilin for several suggestions and fruitful discussions on the presented material and the Research Training Group 1670 "Mathematics Inspired by String Theory and Quantum Field Theory" at the Unversit\"at Hamburg for travel support. Most of the results of this paper have been presented and discussed in talks during the "Research Trimester on Multiple Zeta Values, Multiple Polylogarithms, and Quantum Field Theory" at the ICMAT and therefore the author would also like to thank the organizers of these workshops.  

\section{Bi-brackets and their generating series}
As motivated in the introduction we want to study the following $q$-series:
\begin{dfn} For $r_1,\dots,r_l \geq 0$, $s_1,\dots,s_l > 0$ and we define the following $q$-series
\begin{align*}
 \mb{s_1, \dots , s_l}{r_1,\dots,r_l} :=&\sum_{\substack{u_1 > \dots > u_l> 0 \\ v_1, \dots , v_l >0}} \frac{u_1^{r_1}}{r_1!} \dots \frac{u_l^{r_l}}{r_1!} \cdot \frac{v_1^{s_1-1} \dots v_l^{s_l-1}}{(s_1-1)!\dots(s_l-1)!}   \cdot q^{u_1 v_1 + \dots + u_l v_l} \in \Q[[q]]\,
 \end{align*}
which we call \emph{bi-brackets} of weight $r_1+\dots+r_k+s_1+\dots+s_l$, upper weight $s_1+\dots+s_l$, lower weight $r_1+\dots+r_l$ and length $l$. By $\bMD$ we denote the $\Q$-vector space spanned by all bi-brackets and $1$. 
\end{dfn}

The factorial factors in the definition will become clear when considering their generating functions and the connection to multiple zeta values. For $r_1=\dots=r_l=0$ the bi-brackets are just the brackets 
\[  \mb{s_1, \dots , s_l}{0,\dots,0} = [s_1,\dots,s_l]  \]  
as defined and studied in \cite{BK}. The space spanned by all brackets form a differential $\Q$-algebra $\MD$ with the differential given by $\dif = q \frac{d}{dq}$. We will see that the bi-brackets are also closed under the multiplication of formal power series and therefore $\bMD$ is a $\Q$-algebra with subalgebra $\MD$ (see Theorem \ref{thm:algebra}). 

\begin{dfn}\label{def:genfct}
For the generating function of the bi-brackets we write
\[ \mt{ X_1,\dots,X_l}{Y_1,\dots,Y_l}:= \sum_{\substack{s_1,\dots,s_l > 0 \\ r_1, \dots, r_l > 0}} \mb{s_1\,,\dots\,,s_l}{r_1-1\,,\dots\,,r_l-1} X_1^{s_1-1} \dots X_l^{s_l-1} \cdot Y_1^{r_1-1} \dots Y_l^{r_l-1} \,. \]
These are elements in the ring $\bMDG =\varinjlim_j \bMD[[X_1,\dots,X_j,Y_1,\dots,Y_j]]$ of all generating series of bi-brackets.
\end{dfn}
To derive relations between bi-brackets we will prove functional equations for their generating functions. The key fact for this is that there are two different ways of expressing these given by the following Theorem.
\begin{thm} \label{propgenfct}
For $n\in \N$ set
\[ E_n(X) := e^{nX}\,\quad \text{ and } \quad L_n(X) := \frac{e^X q^n}{1-e^X q^n} \in \Q[[q,X]]\,. \]
Then for all $l\geq 1$ we have the following two different expressions for the generating functions:
\begin{align*}
\mt{X_1,\dots,X_l}{Y_1,\dots,Y_l} &= \sum_{u_1>\dots>u_l>0} \prod_{j=1}^l E_{u_j}(Y_j) L_{u_j}(X_j) \\
&=  \sum_{u_1>\dots>u_l>0} \prod_{j=1}^l E_{u_j}(X_{l+1-j}-X_{l+2-j}) L_{u_j}(Y_1+\dots+Y_{l-j+1})
\end{align*}
(with $X_{l+1} := 0$). In particular the \emph{partition relations} holds:
\begin{equation} \label{eq:partition}
 \mt{X_1,\dots,X_l}{Y_1,\dots,Y_l} \overset{P}{=} \mt{ Y_1 + \dots + Y_l ,\dots ,Y_1 + Y_2, Y_1}{X_l , X_{l-1} - X_l, \dots , X_{1}-X_{2}}\,. 
\end{equation}
\end{thm} 
\begin{prf}
First rewrite the generating function as 
\begin{align*}
  \mt{X_1,\dots,X_l}{Y_1,\dots,Y_l} &=  \sum_{\substack{s_1,\dots,s_l > 0 \\ r_1, \dots, r_l > 0 \\u_1>\dots>u_l > 0 \\ v_1, \dots, v_l > 0}} \prod_{j=1}^l \frac{u_j^{r_j-1}}{(r_j-1)!} \frac{v_j^{s_j-1}}{(s_j-1)!} q^{u_j v_j}   X_j^{s_j-1}  Y_j^{r_j-1} \\
&=  \sum_{\substack{u_1>\dots>u_l > 0 \\ v_1, \dots, v_l > 0}} \prod_{j=1}^l e^{v_j X_j} e^{u_j Y_j} q^{u_j v_j}
 \end{align*}
The first statement follows directly by using the geometric series because
\[ \sum_{v>0} e^{vX} q^{uv} = \frac{e^X q^u}{1-e^X q^u} = L_u(X) \] 
For the second statement set $u_j = u'_1 + \dots + u'_{l-j+1}$ and $v'_j = v_1 + \dots + v_{l-j+1}$ (i.e. $v_j = v'_{l-j+1} - v'_{l-j+2}$ and $v_{l+1}:=0$) for $1 \leq j \leq l$. This gives
\[ q^{u_1v_1 +\dots + u_l v_l} = q^{(u'_1+\dots+u'_l) v_1 + (u'_1 + \dots + u'_{l-1}) v_2 + \dots + u'_1 v_l} =  q^{(v_1+\dots+v_l) u'_1 + \dots + v_1 u'_l }  = q^{v'_1 u'_1 +\dots +  v'_l u'_l}\,  \]
and the summation over $u_1>\dots>u_l > 0$ and $v_1, \dots, v_l > 0$  changes to a summation over $u'_1,\dots,u'_l > 0$ and $v'_1 > \dots > v'_l >0$ and therefore we obtain
\begin{align*}
 \sum_{\substack{u_1>\dots>u_l > 0 \\ v_1, \dots, v_l > 0}} \prod_{j=1}^l e^{v_j X_j} e^{u_j Y_j} q^{u_j v_j} &= \sum_{\substack{v'_1 > \dots > v'_l >0 \\ u'_1,\dots,u'_l > 0}} \prod_{j=1}^l e^{(v'_{l-j+1} - v'_{l-j+2}) X_j} e^{(u'_1 + \dots + u'_{l-j+1}) Y_j} q^{v'_j u'_j} \\
 &=  \sum_{v'_1>\dots>v'_l>0} \prod_{j=1}^l e^{v'_j (X_{l-j+1}-X_{l-j+2})} L_{v'_j}(Y_1+\dots+Y_{l-j+1}) \,
\end{align*}
which is exactly the representation of the generating function. 
\end{prf}
Compare the relation \eqref{eq:partition} to the conjugation $\eqref{eq:conj}$ of partitions given at the end of this section. 

\begin{rem}\begin{enumerate}[i)]
\item The bi-brackets and their generating series also give examples of what is called a bimould by Ecalle in \cite{E}. In his language the equation \eqref{eq:partition} states that the bimould of generating series of bi-brackets is swap invariant. 
\item In \cite{Z} the author studied a variation of the bi-brackets, namely the series
\[ \mathfrak{Z}\mb{s_1,\dots,s_l}{r_1,\dots,r_l} = \sum_{\substack{m_1,\dots,m_l >0 \\ d_1, \dots,d_l >0}} \frac{m_1^{r_1-1} d_1^{s_1-1} \dots m_l^{r_l-1} d_l^{s_l-1} q^{(m_1+\dots+m_l)d_1 + \dots + m_l d_l}}{(r_1-1)!(s_1-1)! \dots (r_l-1)! (s_l-1)!} \,,  \]
which he calls multiple $q$-zeta brackets. These can be written in terms of bi-brackets and vice versa. For this model the equation \eqref{eq:partition}, which in \cite{Z} is called duality, has the nice form 
\[ \mathfrak{Z}\mb{s_1,\dots,s_l}{r_1,\dots,r_l} = \mathfrak{Z}\mb{s_l,\dots,s_1}{r_l,\dots,r_1} \,. \] 
\end{enumerate}
\end{rem}

\begin{cor} (Partition relation in length one and two) 
For $r,r_1,r_2 \geq 0$ and $s,s_1,s_2 > 0$ we have the  following relations in length one and two
\begin{align*}
\mb{s}{r} &= \mb{r+1}{s-1} \,,\\
\mb{s_1,s_2}{r_1,r_2} &= \sum_{\substack{0 \leq j \leq r_1\\0 \leq k \leq s_2-1}} (-1)^k \binom{s_1-1+k}{k} \binom{r_2+j}{j} \mb{r_2+j+1 \,,r_1-j+1}{s_2-1-k \,, s_1-1+k} \,.
\end{align*}
\end{cor}
\begin{prf}
In the smallest cases the Theorem \ref{propgenfct} gives
\begin{align*}
\mt{X}{Y} &= \mt{Y}{X} \quad\text{ and }\quad\mt{X_1,X_2}{Y_1,Y_2} = \mt{Y_1+Y_2,Y_1}{X_2, X_1-X_2} \,.
\end{align*}
The statement follows by considering the coefficients of $X^{s-1} Y^{r}$ and $X_1^{s_1-1} X_2^{s_2-1} Y^{r_1} Y^{r_2}$ in these equations.
\end{prf}

\begin{ex}
\begin{enumerate}[i)]
\item Some examples for the length two case: 
\begin{align*}
\mb{1,1}{1,1} &= \mb{2,2}{0,0} + 2\mb{3,1}{0,0} \,,\quad \mb{3,3}{0,0} = 6 \mb{1,1}{0,4}-3 \mb{1,1}{1,3}+\mb{1,1}{2,2} \,, \\
\mb{2,2}{1,1} &=-2 \mb{2,2}{0,2}+\mb{2,2}{1,1}-4 \mb{3,1}{0,2}+2 \mb{3,1}{1,1} \,,\\ 
\mb{1,2}{2,3} &=-\mb{4,3}{0,1}+\mb{4,3}{1,0}-4 \mb{5,2}{0,1}+4 \mb{5,2}{1,0}-10 \mb{6,1}{0,1}+10 \mb{6,1}{1,0} \,.
\end{align*}
\item Another family of relations which can be obtained by the partition relation is
\[ \mb{\{ 1 \}^n }{\{0\}^{j-1},1,\{0\}^{n-j}} = \sum_{k=1}^{n-j+1} [\{1\}^{k-1} , 2, \{1\}^{n-k}] \]
for $1 \leq j \leq n$. For example: 
\[ \mb{1}{1} = [2]\,, \quad \mb{1,1,1}{0,1,0} = [1,2,1] + [2,1,1] \,. \]
\end{enumerate}
\end{ex}

\begin{rem}
We end the discussion on bi-brackets and their generating series by interpreting the coefficients of the bi-brackets as weighted sums over partitions which gives an natural explanation for the partition relation \eqref{eq:partition}. By a partition of a natural number $n$ with $l$ parts we denote a representation of $n$ as a sum of $l$ distinct natural numbers, i.e. $15 = 4 + 4 + 3 + 2 + 1 + 1$ is a partition of $15$ with the $4$ parts given by $4,3,2,1$. We identify such a partition with a tuple $(u,v) \in \N^l \times \N^l$ where the $u_j$'s are the $l$ distinct numbers in the partition and the $v_j$'s count their appearance in the sum. The above partition of $15$ is therefore given by the tuple $(u,v) = ((4,3,2,1),(2,1,1,2))$. By $P_l(n)$ we denote all partitions of $n$ with $l$ parts and hence we set
\[ P_l(n) := \left\{ (u,v) \in \N^l \times \N^l \, \mid \, n = u_1 v_1 + \dots + u_l v_l \, \text{ and } \, u_1 > \dots > u_l > 0  \right\}   \] 
On the set $P_l(n)$ one has an involution given by the conjugation $\rho$ of partitions which can be obtained by reflecting the corresponding Young diagram across the main diagonal. 
\begin{figure}[h!]
\captionsetup{width=0.8\textwidth}
\[ {\tiny ((4,3,2,1),(2,1,1,2)) = \yng(4,4,3,2,1,1)  \quad  \overset{\rho}{\xrightarrow{\hspace*{1cm}}} \quad \yng(6,4,3,2) = ((6,4,3,2),(1,1,1,1))} \]
  \caption{The conjugation of the partition $15 = 4 + 4 + 3 + 2 + 1 + 1$ is given by $\rho( ((4,3,2,1),(2,1,1,2)) ) =  ((6,4,3,2),(1,1,1,1))$ which can be seen by reflection the corresponding Young diagram at the main diagonal.  }
\end{figure}
On the set $P_l(n)$ the conjugation $\rho$ is explicitly given by $\rho( (u,v) ) = (u',v')$ where $u'_j = v_1 + \dots + v_{l-j+1}$ and $v'_j = u_{l-j+1}-u_{l-j+2}$ with $u_{l+1} := 0$, i.e.  
\begin{equation} \label{eq:conj}
\rho : \binom{u_1, \dots , u_l}{v_1, \dots ,v_l} {\longmapsto} \binom{v_1+\dots+v_l, \dots,v_1+v_2,v_1}{u_l,u_{l-1}-u_l,\dots,u_1-u_2} \,.
\end{equation}
By the definition of the bi-brackets its clear that with the above notation they can be written as
\begin{align*}
 \mb{s_1, \dots , s_l}{r_1,\dots,r_l} :=& \frac{1}{r_1! (s_1-1)! \dots r_l! (s_l-1)!} \sum_{n>0} \left(  \sum_{(u,v) \in P_l(n)} u_1^{r_1} v_1^{s_1-1} \dots  u_l^{r_l} v_l^{s_l-1}  \right) q^n \,. 
 \end{align*}
The coefficients are given by a sum over all elements in $P_l(n)$ and therefore it is invariant under the action of $\rho$. As an example consider $[2,2]$ and apply $\rho$ to the sum then we obtain 
\begin{align} \label{eq:2200ex}
\begin{split}
[2,2] &=  \sum_{n>0} \left( \sum_{(u,v) \in P_2(n)}  v_1 \cdot v_2 \right) q^n = \sum_{n>0} \left(\sum_{ \rho((u,v))=(u',v') \in P_2(n)}  u_2' \cdot (u_1'-u_2') \right) q^n \\
&= \sum_{n>0} \left( \sum_{(u',v') \in P_2(n)}  u_2' \cdot u_1' \right) q^n - \sum_{n>0} \left( \sum_{(u',v') \in P_2(n)}  u_2'^2 \right) q^n = \mb{1,1}{1,1} - 2 \mb{1,1}{0,2} \,.
\end{split}
\end{align}
This is exactly the relation one obtains by using the partition relation. Another trivial connection to partitions is given by the following: The coefficients of the brackets of the form $[\{1\}^l]$ count the number of partitions of length $l$. Summing over all length one therefore obtains the generating functions of all partitions: 
\[ \sum_{l>0} [\{1\}^l] = \sum_{n>0} p(n) q^n =  \prod_{k=1}^\infty \frac{1}{1-q^k} \,.\]
\end{rem}

\section{The algebra of bi-brackets}

The partition relations give relations in a fixed length. To obtain relations with mixed length we need to consider the algebra structure on the space $\bMD$. For this we first consider the product of bi-brackets in length one and then use the algebraic setup of quasi-shuffle algebras for the arbitrary length case. 
\begin{lem} \label{lem:lnprod}
Let $B_k$ be the $k$-th Bernoulli number, then we get for all $n\in \N$ 
\begin{align*}
 &L_n(X)\cdot L_n(Y) =\sum_{k>0} \frac{B_k}{k!}(X-Y)^{k-1} L_n(X) + \sum_{k>0} \frac{B_k}{k!}(Y-X)^{k-1} L_n(Y) + \frac{L_n(X) - L_n(Y)}{X-Y} \,. 
\end{align*}
\end{lem}
\begin{prf}
By direct computations one obtains
\[ L(X) \cdot L(Y) = \frac{1}{e^{X-Y}-1} L(X) + \frac{1}{e^{Y-X}-1} L(Y) \,. \]
The statement follows  then by the definition of the Bernoulli numbers 
\[ \frac{X}{e^X-1} = \sum_{n\geq0} \frac{B_n}{n!} X^n \,.\]
\end{prf}

\begin{lem} \label{lem:genprod11}
The product of two generating functions in length one can be written as 
\begin{enumerate}[i)]
\item ("Stuffle product for bi-brackets")
\begin{align*}
\mt{X_1}{Y_1} \cdot \mt{X_2}{Y_2}&=
\mt{X_1,X_2}{Y_1,Y_2} + \mt{X_2,X_1}{Y_2,Y_1} + \frac{1}{X_1-X_2} \left( \mt{X_1}{Y_1+Y_2} - \mt{X_2}{Y_1+Y_2} \right) \\
&+ \sum_{k=1}^\infty \frac{B_k}{k!} (X_1-X_2)^{k-1} \left( \mt{X_1}{Y_1+Y_2} + (-1)^{k-1} \mt{X_2}{Y_1+Y_2} \right)  \,. 
\end{align*}
\item ("Shuffle product for bi-brackets") 
\begin{align*}
\mt{X_1}{Y_1} \cdot \mt{X_2}{Y_2}&= \mt{X_1+X_2, X_1}{Y_2 \,, Y_1 - Y_2} + \mt{X_1+X_2, X_2}{Y_1 , Y_2 - Y_1} + \frac{1}{Y_1-Y_2} \left( \mt{X_1+X_2}{Y_1} - \mt{X_1+X_2}{Y_2} \right) \\
&+ \sum_{k=1}^\infty \frac{B_k}{k!} (Y_1-Y_2)^{k-1} \left( \mt{X_1+X_2}{Y_1} + (-1)^{k-1} \mt{X_1+X_2}{Y_2} \right)  
\end{align*}
\end{enumerate}
\end{lem}
\begin{prf}
We prove $i)$ and $ii)$ by using the two different ways of writing the generating functions given by Theorem \ref{propgenfct}.
\begin{enumerate}[i)]
\item
By direct calculation it is
\begin{align*}
\mt{X_1}{Y_1}\cdot &\mt{X_2}{Y_2}=\mt{X_1,X_2}{Y_1,Y_2} + \mt{X_2,X_1}{Y_2,Y_1} + \sum_{n>0} E_n(Y_1+Y_2) L_n(X_1) L_n(X_2)\,. 
\end{align*}
Applying the \lemref{lem:lnprod} to the last term yields the statement.
\item 
The partition relation in length one and two ($P$) in \eqref{eq:partition} states
\begin{align*}
\mt{X_1}{Y_1} \overset{P}{=} \mt{Y_1}{X_1} \,,\quad \mt{X_1,X_2}{Y_1,Y_2} \overset{P}{=} \mt{Y_1+Y_2,Y_1}{X_2,X_1-X_2} \,,
\end{align*}
and together with  i) we obtain
\begin{align*} 
\mt{X_1}{Y_1} \cdot \mt{X_2}{Y_2}& \overset{P}{=}\mt{Y_1}{X_1} \cdot \mt{Y_2}{X_2} \overset{i)}{=} \mt{Y_1,Y_2}{X_1,X_2} + \mt{Y_2,Y_1}{X_2,X_1} + \frac{1}{Y_1-Y_2} \left( \mt{Y_1}{X_1+X_2} - \mt{Y_2}{X_1+X_2} \right) \\
&+ \sum_{k=1}^\infty \frac{B_k}{k!} (Y_1-Y_2)^{k-1} \left( \mt{Y_1}{X_1+X_2} + (-1)^{k-1} \mt{Y_2}{X_1+X_2} \right) \\
&\overset{P}{=}\mt{X_1+X_2, X_1}{Y_2 \,, Y_1 - Y_2} + \mt{X_1+X_2, X_2}{Y_1 , Y_2 - Y_1} + \frac{1}{Y_1-Y_2} \left( \mt{X_1+X_2}{Y_1} - \mt{X_1+X_2}{Y_2} \right) \\
&+ \sum_{k=1}^\infty \frac{B_k}{k!} (Y_1-Y_2)^{k-1} \left( \mt{X_1+X_2}{Y_1} + (-1)^{k-1} \mt{X_1+X_2}{Y_2} \right)  \,.
\end{align*}
\end{enumerate}
\end{prf}

\begin{prop}\label{prop:expstsh} 
For $s_1,s_2 > 0$ and $r_1,r_2 \geq 0$ we have the following two expressions for the product of two bi-brackets of length one:
\begin{enumerate}[i)]
\item ("Stuffle product for bi-brackets")
\begin{align*}
\mb{s_1}{r_1} \cdot \mb{s_2}{r_2} &=\mb{s_1,s_2}{r_1,r_2} + \mb{s_2,s_1}{r_2,r_1} + \binom{r_1+r_2}{r_1}\mb{s_1+s_2}{r_1+r_2} \\
&+\binom{r_1+r_2}{r_1}\sum_{j=1}^{s_1} \frac{(-1)^{s_2-1} B_{s_1+s_2-j}}{(s_1+s_2-j)!} \binom{s_1+s_2-j-1}{s_1-j}  \mb{j}{r_1+r_2} \\
&+\binom{r_1+r_2}{r_1}  \sum_{j=1}^{s_2} \frac{(-1)^{s_1-1} B_{s_1+s_2-j}}{(s_1+s_2-j)!} \binom{s_1+s_2-j-1}{s_2-j} \mb{j}{r_1+r_2} 
\end{align*}

\item ("Shuffle product for bi-brackets") 
\begin{align*}
\mb{s_1}{r_1} \cdot \mb{s_2}{r_2} &= \sum_{\substack{1 \leq j \leq s_1 \\ 0 \leq k \leq r_2}} \binom{s_1+s_2-j-1}{s_1-j} \binom{r_1+r_2-k}{r_1} (-1)^{r_2-k} \mb{s_1+s_2-j, j}{k,r_1+r_2-k} \\
&+\sum_{\substack{1\leq j \leq s_2 \\ 0 \leq k \leq r_1}} \binom{s_1+s_2-j-1}{s_1-1} \binom{r_1+r_2-k}{r_1-k} (-1)^{r_1-k} \mb{s_1+s_2-j, j}{k,r_1+r_2-k}\\
&+\binom{s_1+s_2-2}{s_1-1} \mb{s_1+s_2-1}{r_1+r_2+1} \\
&+\binom{s_1+s_2-2}{s_1-1}\sum_{j=0}^{r_1} \frac{(-1)^{r_2} B_{r_1+r_2-j+1}}{(r_1+r_2-j+1)!} \binom{r_1+r_2-j}{r_1-j}  \mb{s_1+s_2-1}{j} \\
&+\binom{s_1+s_2-2}{s_1-1}\sum_{j=0}^{r_2} \frac{(-1)^{r_1} B_{r_1+r_2-j+1}}{(r_1+r_2-j+1)!} \binom{r_1+r_2-j}{r_2-j}  \mb{s_1+s_2-1}{j} 
\end{align*}

\end{enumerate}
\end{prop}

\begin{prf}

\begin{enumerate}[i)]
\item 
By \lemref{lem:genprod11} it is

\begin{align*}
\mt{X_1}{Y_1} \cdot \mt{X_2}{Y_2}&=\underbrace{\mt{X_1,X_2}{Y_1,Y_2} + \mt{X_2,X_1}{Y_2,Y_1}}_{=: T_1} +\underbrace{ \frac{1}{X_1-X_2} \left( \mt{X_1}{Y_1+Y_2} - \mt{X_2}{Y_1+Y_2} \right)}_{=: T_2} \\
&+\underbrace{ \sum_{k=1}^\infty \frac{B_k}{k!} (X_1-X_2)^{k-1} \left( \mt{X_1}{Y_1+Y_2} + (-1)^{k-1} \mt{X_2}{Y_1+Y_2} \right)  }_{=: T_3}\,. 
\end{align*}

We are going to calculate the coefficient of $X^{s_1-} X_2^{s_2-1} Y_1^{r_1} Y_2^{r_2}$ in this equation. Clearly $\mb{s_1,s_2}{r_1,r_2} + \mb{s_2,s_1}{r_2,r_1}$ is the coefficient of $T_1$ and by the use of
\[ \sum_{s>0} c_s \frac{X_1^{s-1} - X_2^{s-1}}{X_1-X_2} =   \sum_{s>0} c_s \sum_{j=0}^{s-2} X_1^{s-2-j} X_2^j  = \sum_{a,b >0} c_{a+b} X_1^{a-1} X_2^{b-1} \]
one obtains
\begin{align*}
T_2&=\frac{1}{X_1-X_2} \left( \mt{X_1}{Y_1+Y_2} - \mt{X_2}{Y_1+Y_2} \right) = \sum_{s_1,s_2,r > 0} \mb{s_1+s_2}{r-1} X_1^{s_1-1} X_2^{s_2-1} (Y_1+Y_2)^{r-1} \\
&= \sum_{\substack{s_1,s_2 > 0\\r_1,r_2 > 0}} \binom{r_1+r_2-2}{r_1-1}\mb{s_1+s_2}{r_1+r_2-2} X_1^{s_1-1} X_2^{s_2-1} Y_1^{r_1-1} Y_2^{r_2-1}  \,.
\end{align*}
With a bit more tedious but similar calculation one shows that the remaining terms are the coefficients of $T_3$.
\item  This statement follows by a similar calculation as in i). 
\end{enumerate}
\end{prf}

We now want to recall the algebraic setting of Hoffman for quasi-shuffle products and give the necessary notations for the rest of the paper. 
\begin{dfn} 
Let $A$ (the alphabet) be a countable set of letters, $\Q A$ the $\Q$-vector space generated by these letters and $\Q\langle A \rangle$ the noncommutative polynomial algebra over $\Q$ generated by words with letters in $A$. For a commutative and associative product $\diamond$ on $\Q A$, $a,b \in A$ and $w,v \in \Q\langle A\rangle$ we define on $\Q\langle A\rangle$ recursively a product by $1 \odot w=w \odot 1=w$ and
\begin{equation} \label{qshp} aw \odot bv := a(w \odot bv) + b(aw \odot v) + (a \diamond b)(w \odot v) \,. \end{equation}
 By a result of Hoffman (\cite{H1}) $(\Q\langle A \rangle,\odot)$ is a commutative $\Q$-algebra which is called a \emph{quasi-shuffle algebra}.
\end{dfn}

\begin{nota}
Let us now recall some basic notations for the shuffle and the stuffle product which are the easiest examples of quasi-shuffle products. Since we will deal with the shuffle product for different alphabets simultaneously we will use some additional notations for this. For the alphabet $A_{xy} := \{ x,y \}$ set $\h = \Q\langle A_{xy} \rangle$ and $\h^1 = \h y$. It is easy to see that $\h^1$ is generated by the elements $z_j = x^{j-1} y$ with $j \in \N$, i.e. $\h^1 = \Q\langle A_{z} \rangle$ with $A_z := \{ z_1,z_2, \dots \}$. By $|w|$ we denote the the weight of a word $w\in \h$ which is given by the number of letters (in the alphabet $A_{xy}$)  of $w$. On $\h^1$ we have the following two products with respect to the alphabet $A_z$ which we call the \emph{index-shuffle}, denoted by $\shuffle$ with $\diamond \equiv 0$, and the \emph{stuffle product}, denoted by $\ast$ with $z_j \diamond z_i = z_{j+i}$, i.e. we have for $a,b \in \N$ and  $w,v \in \h^1$:
\begin{align}
\begin{split}
 z_a w \sh z_b v &= z_a(w \sh z_b v) + z_b(z_a w \sh v)  \,, \\
 z_a w \ast z_b v &= z_a(w \ast z_b v) + z_b(z_a w \ast v)  + z_{a+b}( w \ast v) \,.
\end{split}
\end{align}
By $(\h_z^1, \sh)$  and $(\h^1_z, \ast)$ we denote the corresponding $\Q$-algebras, where the subscript $z$ indicates that we consider the quasi-shuffle with respect to the alphabet $A_z$. We can also define the shuffle product on $\h^1$ with respect to the alphabet $A_{xy}$, which we call the \emph{shuffle product}, and by $(\h_{xy}^1,\sh)$ we denote the corresponding $\Q$-algebra. 
\end{nota}

We now want to find a $\diamond$ and a suitable alphabet such that we can view the algebra of bi-brackets as a quasi-shuffle algebra.
For $a,b \in \N$ define the numbers $\lambda^j_{a,b}  \in \Q$ for $1 \leq j \leq a$ as
\[ \lambda^j_{a,b} = (-1)^{b-1} \binom{a+b-j-1}{a-j} \frac{B_{a+b-j}}{(a+b-j)!} \,. \]
For the alphabet $A^{\text{bi}}_{z} := \{z_{s,r} \mid s,r\in \Z \,, s\geq 1\,, r\geq 0\}$ we define on $\Q A^{\text{bi}}_{z}$ the product
\begin{align*}
 z_{s_1,r_1} \boxcircle z_{s_2,r_2} =&  \binom{r_1+r_2}{r_1} \sum_{j=1}^{s_1} \lambda^j_{a,b} z_{j,r_1+r_2} + \binom{r_1+r_2}{r_1}\sum_{j=1}^{s_2} \lambda^j_{b,a}  z_{j,r_1+r_2}  \\\
 &+\binom{r_1+r_2}{r_1} z_{s_1+s_2,r_1+r_2} 
\end{align*} 
and on $\Q\langle A^{\text{bi}}_{z} \rangle$ the quasi-shuffle product
 \begin{align*}
  z_{s_1,r_1} w \boxast z_{s_2,r_2} v &= z_{s_1,r_1}(w \boxast z_{s_2,r_2} v) + z_{s_2,r_2}(z_{s_1,r_1} w \boxast v)  +  (z_{s_1,r_1} \boxcircle z_{s_2,r_2})(w \boxast v) \,. \\
\end{align*}

\begin{thm}\label{thm:algebra}
\begin{enumerate}[i)]
\item The product $\boxcircle$ on $\Q A^{\text{bi}}_{z}$ is associative and therefore  $(\Q\langle A^{\text{bi}}_{z}\rangle,\boxast)$ is a quasi-shuffle Algebra. 
\item The map $\mb{.}{}: (\Q\langle A^{\text{bi}}_{z}\rangle,\boxast) \rightarrow \left( \bMD , \cdot \right)$ given by
\[  w = z_{s_1,r_1} \dots  z_{s_l,r_l} \longmapsto [w] = \mb{s_1,\dots,s_l}{r_1,\dots,r_l} \]
fulfills $[w \boxast v]=[w] \cdot [v]$ and therefore $\bMD$ is a $\Q$-algebra.
\end{enumerate}
\end{thm}
\begin{prf}
 Using Proposition 2.3 in \cite{BK} it is easy to see that
\begin{equation}\label{eq:polylogdst}
 \mb{s_1,\dots,s_l}{r_1, \dots , r_l} = \sum_{u_1 > \dots > u_l > 0} \frac{u_1^{r_1}}{r_1!}  \Lit_{s_1}(q^{u_1}) \dots  \frac{u_l^{r_l}}{r_1!} \Lit_{s_l}(q^{u_l})\,, 
 \end{equation}
where $\Lit_s(x) = \frac{ \Li_{1-s}(x)}{(s-1)!}$. Due to Lemma \ref{lem:lnprod} (see also Lemma 2.5 in \cite{BK}) we have
\[ \Lit_{a}(z) \cdot  \Lit_{b}(z) = \sum_{j=1}^a \lambda^j_{a,b} \Lit_{j}(z) + \sum_{j=1}^b \lambda^j_{b,a}  \Lit_{j}(z)  + \Lit_{a+b}(z) \,, \]
This proves the first statement and the second statement follows directly by the definition of $\boxast$.  
\end{prf}

\begin{rem} \label{rmk:pstp}
As we saw in the proof of Proposition \ref{lem:genprod11} for the product of two length one bi-brackets, the shuffle product of bi-brackets is obtained by applying the partition relation, the stuffle product and again the partition relation. This of course works for arbitrary lengths and yields a natural way to obtain the shuffle product for bi-brackets. To make this precise denote by $P:  \Q\langle A^{\text{bi}}_{z}\rangle \rightarrow  \Q\langle A^{\text{bi}}_{z}\rangle$ the linearly extended map which sends a word $w = z_{s_1,r_1} \dots z_{s_l,r_l}$ to the linear combination of words corresponding to the partition relation. Using this convention the shuffle product for brackets can be written in $ \Q\langle A^{\text{bi}}_{z}\rangle$ for two words $u , v \in \Q\langle A^{\text{bi}}_{z}\rangle$ as $ P\left( P(u) \boxast P(v) \right)$, i.e. the stuffle and shuffle product for bi-brackets can be written as
\begin{equation}\label{eq:stsh}
 [u] \cdot [v] \overset{st}{=} [ u \boxast v] \,, \qquad [u] \cdot [v] \overset{sh}{=} [ P\left( P(u) \boxast P(v) \right) ] \,.
\end{equation}
\end{rem}

\begin{rem}
As mentioned in the introduction the bi-brackets can be seen as a $q$-analogue of MZV: Define for $k\in \N$ the map $\Q[[q]] \rightarrow \R \cup \{ \infty \}$ by $Z_k(f) = \lim_{q \to 1} (1-q)^{k}f(q)$, which was introduced and discussed in \cite{BK} for the subspace $\MD \subset \Q[[q]]$. On the bi-brackets this map is given by the following: Assume that $s_1 > r_1+1$ and $s_j \geq r_j + 1$ for $j=2, \dots, l$, then, using the description \eqref{eq:polylogdst} (see eg. Proposition 1 in \cite{Z}), we obtain
\[ Z_{s_1+\dots+s_l}\mb{s_1,\dots,s_l}{r_1,\dots,r_l} = \frac{1}{r_1! \dots r_l!} \zeta(s_1-r_1,\dots,s_l-r_l) \,.\] 
Even though we don't want to discuss this issue in this note it is worth mentioning that an other motivation for considering the bi-brackets was to describe the kernel of the map $Z_k$ on the space $\grw_k \MD$. This connection will be subject of upcoming works. Applying the map $Z_k$ to the equation \eqref{eq:stsh} one obtains the stuffle and shuffle product formula for MZV (See \cite{Z}).  Finally we want mention that there are several other different types of $q$-analogues which also have a stuffle and shuffle like structure (See for example \cite{JMK} and \cite{Zh} for a nice overview). 
\end{rem}

\section{Derivatives and modular forms}
In this section we want to discuss derivatives of bi-brackets with respect to the differential operator $q \frac{d}{dq}$ and their connections to modular forms. For this we first introduce the following notations: 
\begin{dfn}
On $\bMD$ we have the increasing filtrations $\filw_{\bullet}$
given by the upper weight,$\fild_{\bullet}$ give by the lower weight and $\fille_{\bullet}$
given by the length, i.e., we have for $A \subseteq \bMD$
\begin{align*}
\filw_k(A) &:=  \big< \mb{s_1,\dots,s_l}{r_1, \dots, r_l} \in A \,\big|\, 0 \leq l \leq k \,,\, s_1+\dots+s_l \le k \,\big>_{\Q}\\
\fild_k(A) &:=  \big< \mb{s_1,\dots,s_l}{r_1, \dots, r_l} \in A \,\big|\, 0 \leq l \leq k \,,\, r_1+\dots+r_l \le k \,\big>_{\Q}\\
\fille_l(A) &:=  \big<\mb{s_1,\dots,s_t}{r_1, \dots, r_t} \in A \,\big|\, t\le l \,\big>_{\Q}\,.
\end{align*}
If we consider the length and weight filtration at the same time we use the short notation $\filwle_{k,l} := \filw_k \fille_l$ and similar for the other filtrations.
\end{dfn}

\begin{prop}\label{prop:derivative}
Let $\dif := q \frac{d}{dq}$ then we have
\[ \dif \mb{s_1, \dots , s_l}{r_1,\dots,r_l} = \sum_{j=1}^l \left( s_j (r_j+1) \mb{s_1\,,\dots\,,s_{j-1}\,,s_j+1\,,s_{j+1},\dots \,,s_l}{r_1 \,, \dots \,,r_{j-1}\,,r_j+1\,,r_{j+1} \,,\dots \,, r_l} \right) \, \]
and therefore $\dif \left( \operatorname{Fil}^{\operatorname{W},\operatorname{D},\operatorname{L}}_{k,d,l}(\bMD) \right) \subset  \operatorname{Fil}^{\operatorname{W},\operatorname{D},\operatorname{L}}_{k+1,d+1,l}(\bMD) $.
\end{prop}
\begin{prf}
This is an easy consequence of the definition of bi-brackets and the fact that $\dif \sum_{n>0} a_n q^n = \sum_{n>0} n a_n q^n$. Another way to see this is by the fact that the operator $\dif$ on the generating series of bi-brackets can be written as 
\[ \dif \mt{ X_1,\dots,X_l}{Y_1,\dots,Y_l} = \sum_{j=1}^l   \frac{\partial}{\partial X_j} \frac{\partial}{\partial Y_j}  \mt{ X_1,\dots,X_l}{Y_1,\dots,Y_l}\,,  \] 
which follows from 
\[ \dif E_n(Y) L_n(X) = \dif \frac{e^{nY} e^{X} q^n}{(1-e^X q^n)} =   \frac{n e^{nY} e^{X} q^n}{(1-e^X q^n)^2} = \frac{\partial}{\partial X} \frac{\partial}{\partial Y} E_n(Y) L_n(X)\,.\] 
\end{prf}
Proposition \ref{prop:derivative} suggests that the bi-brackets can be somehow viewed as  partial derivatives of the brackets with total differential $\dif$. In this part we want to give some explicit results on the following conjecture which was inspired by numerical experiments and which, with the above interpretation, states that the space $\MD$ is closed under partial derivatives.
\begin{conj}\label{conj:bdmd}
The algebra $\bMD$ of bi-brackets is a subalgebra of $\MD$ and in particular we have
\[ \operatorname{Fil}^{\operatorname{W},\operatorname{D},\operatorname{L}}_{k,d,l}(\bMD) \subset  \operatorname{Fil}^{\operatorname{W},\operatorname{L}}_{k+d,l+d}(\MD) \,. \]
\end{conj}
\begin{prop}\label{prop:conjinlen1} For $l=1$ the conjecture \ref{conj:bdmd} is true. 
\end{prop}
\begin{prf} In \cite{BK} the authors proved that $\dif \MD \subset \MD$. Due to Proposition \ref{prop:derivative} we therefore have $\mb{s}{r} \in \MD$, i.e. the Conjecture is true for the length one case. 
\end{prf}

\begin{rem}
In $\cite{BK2}$ it will be shown that up to weight $k\leq 7$ every bi-bracket can be written in terms of brackets, by giving upper bounds for the number of algebra generators of bi-brackets. 
\end{rem} 

For lower weight $d=1$ Proposition \ref{prop:conjinlen1} is given explicitly by the following reformulation of Proposition 3.3 in \cite{BK}. 
\begin{prop}
For all $k \geq 1$ it is
\begin{align*}
\mb{k}{1} &= [k] \cdot [1] -\sum_{a+b=k+1} [a,b] - [k,1] + [k] \\
&= [k+1] + \frac{1}{2} [k] - \sum_{\substack{a+b=k+1 \\ a>1}} [a,b] + \sum_{j=2}^{k-1} \frac{B_{k-j+1}}{(k-j+1)!} [j] - \frac{1}{2}\delta_{k,1} [1]  \in \operatorname{Fil}^{\operatorname{W},\operatorname{L}}_{k+1,2}(\MD)
\end{align*}
\end{prop}

\begin{prf}
The functions $L_n(X)$ in the generating function fullfil the following differential equation.
\[ \frac{\partial}{\partial X} L_n(X) = L_n(X)^2 + L_n(X)\,. \]
Therefore we get
\begin{align*}
\frac{\partial}{\partial Y} \mt{X}{Y} =  \sum_{n>0} e^{nX} L_n(Y)^2+ \sum_{n>0} e^{nX} L_n(Y) =  \sum_{n>0} e^{nX} L_n(Y)^2+\mt{X}{Y} \,. 
\end{align*}
The first term also appears in the product of two generating functions:
\begin{align*}
\mt{X}{Y} \cdot \mt{0}{Y} &= \sum_{n_1 > n_2 > 0} e^{n_1 X} L_{n_1}(Y) L_{n_2}(Y) +  \sum_{n_2 > n_1 > 0} e^{n_1 X} L_{n_1}(Y) L_{n_2}(Y) +  \sum_{n>0} e^{nX} L_n(Y)^2 \\
&= \mt{Y,Y}{X,0} + \mt{Y,Y}{0,X} +   \sum_{n>0} e^{nX} L_n(Y)^2 = \mt{X,X}{Y,0} + \mt{X,0}{Y,0} +   \sum_{n>0} e^{nX} L_n(Y)^2 \,.
\end{align*}
And therefore we obtain 
\begin{equation}\label{eq:dylen1} \frac{\partial}{\partial Y} \mt{X}{Y} =\mt{X}{Y} \cdot \mt{0}{Y}-  \mt{X,X}{Y,0} - \mt{X,0}{Y,0} +\mt{X}{Y}\,, \end{equation}
which gives the first expression by considering the coefficient of $X^{k-1}$ in this equation. 
The second statement follows from the explicit stuffle product for bi-brackets in Proposition \ref{prop:expstsh}:
\[ [k] \cdot [1] = [k,1]+[1,k] + [k+1] + \sum_{j=2}^k \frac{B_{k-j+1}}{(k-j+1)!} [j] - \delta_{k,1} [1] \,. \]
\end{prf}

There it not much known so far for the length two and arbitrary weight case of the Conjecture \ref{conj:bdmd}. Using the shuffle brackets we will prove (see Proposition \ref{prop:len2conj} ) that for all $s_1,s_2 \geq 1$ it is
\begin{align*}
\mb{s_1,s_2}{1,0}, \mb{s_1,s_2}{0,1} \in \operatorname{Fil}^{\operatorname{W},\operatorname{L}}_{s_1+s_2+1,3}(\MD) 
\end{align*}
It would be interesting to know whether the approach in the proof of proposition \ref{prop:len2conj} also works for higher lengths, or higher lower weight.


One motivation of considering (bi-)brackets is to build a connection between multiple zeta values and modular forms. In the following we will show how to use the double shuffle structure on the space of bi-brackets described above to prove relations between modular forms. On the other hand we use results of modular forms to prove relations between bi-brackets. For $k\in \N$ denote by  
\[ \widetilde{G}_{k} =  \frac{\zeta(k)}{(-2 \pi i)^k} + \frac{1}{(k-1)!} \sum_{n>0} \sigma_{k-1}(n) q^n  =     \frac{\zeta(k)}{(-2 \pi i)^k}   + [k] \,. \]
the Eisenstein series of weight $k$. For even $k=2n$ due to Euler we have $\zeta(2n) = \frac{(-1)^{n-1} B_{2n} (2\pi)^{2n}}{2(2n)!}$ and therefore $\widetilde{G}_{2n} =   -\frac{1}{2} \frac{B_{2n}}{(2n)!}  + [2n] =: \beta_{2n} + [2n] \in \filw_{2n}(\MD)$,
for example
\[ \widetilde{G}_2 = -\frac{1}{24} + [2] \,,\quad \widetilde{G}_4 = \frac{1}{1440} + [4] \,, \quad \widetilde{G}_6 = -\frac{1}{60480} + [6] \,.\]

\begin{prop} \label{prop:mf-neu} 
\begin{enumerate}[i)]
\item The ring of modular forms $M(\Gamma_1)$ for $\Gamma_1 = SL_2(\Z)$ 
and the ring of quasi-modular forms 
$\widetilde{M}(\Gamma_1)$ are graded subalgebras of $\MD$.
\item The $\Q$-algebra of quasi-modular forms   
$\widetilde M_k(\Gamma_1)$ 
is closed under the derivation $\dif$ and therefore it is a subalgebra of the graded differential algebra 
$(\MD, \dif)$. 

\item We have the following inclusions of $\Q$-algebras
\[ M_k(\Gamma_1) \subset  \widetilde{M}(\Gamma_1) \subset \MDA \subset \MD \subset \bMD \,. \]  
\end{enumerate}
\end{prop}
\begin{prf} Let $M_k(\Gamma_1)$ (resp.  $\widetilde M_k(\Gamma_1)$) be the space of (quasi-)modular forms of weight $k$ for $\Gamma_1$. Then the first claim follows directly from the well-known facts $M(\Gamma_1) = \bigoplus_{k>1} M(\Gamma_1)_k = \Q[\widetilde{G}_4,\widetilde{G}_6]$ and $\widetilde{M}(\Gamma_1) =\bigoplus_{k>1} \widetilde{M}(\Gamma_1)_k = \Q[\widetilde{G}_2,\widetilde{G}_4,\widetilde{G}_6]$. The second claim is a well known fact in the theory of quasi-modular forms and a proof can be found in \cite{dz} p. 49. It suffices to show that the
 derivatives of the generators are given by
\begin{align*} 
\dif \widetilde{G}_2 &=\dif [2] = 5 \widetilde{G}_4 - 2 \widetilde{G}_2^2 \,,\qquad \dif \widetilde{G}_4 = 15 \widetilde{G}_6 - 8 \widetilde{G}_2 \widetilde{G}_4 \,,  \\
\dif \widetilde{G}_6 &= 20 \widetilde{G}_8 - 12 \widetilde{G}_2 \widetilde{G}_6 = \frac{120}{7} \widetilde{G}_4^2 - 12 \widetilde{G}_2 \widetilde{G}_6 \,,
\end{align*} 
which can be easily shown by the double shuffle relations of bi-brackets. 
\end{prf}

It is a well-known fact from the theory of modular forms that $\widetilde{G}_4^2 = \frac{7}{6} \widetilde{G}_8$ because the space of weight $8$ modular forms for $\Sl_2(\Z)$ is one dimensional. We therefore have   
\begin{align*}
   \frac{1}{720} [4]  + [4] \cdot [4] 
&= \frac{7}{6} [8]\,.
\end{align*}
Using the explicit stuffle product we get
\begin{align*}
 [4] \cdot [4] 
&= 2 [4,4] + [8] + \frac{1}{360} [4] - \frac{1}{1512} [2] \,,
\end{align*}
which then gives the following relation in $\filw_{8}(\MD)$:
\begin{align}\label{eq:relwt8}
[8] = \frac{1}{40} [4] - \frac{1}{252} [2] + 12 [4,4]\,.  
\end{align}

The identity \eqref{eq:relwt8} can also be proven by using the double shuffle relations, i.e. $\widetilde{G}_4^2 = \frac{7}{6} G_8$ can be proven since it is equivalent to it. One can check that 
\begin{align*}
\frac{1}{40} [4] - \frac{1}{252} [2] + 12 [4,4] - [8] = -4 ( [3] \overset{st}{\cdot} [5] -  [3] \overset{sh}{\cdot} [5] ) + 3 ([4] \overset{st}{\cdot} [4] -  [4] \overset{sh}{\cdot} [4])  \,,
\end{align*}
where the right hand side is clearly zero. This purely combinatorial approach to prove this kind of relation is similar to the one in \cite{S}.

Let us now use the theory of modular forms to obtain relations between bi-brackets. 
It is a well-known fact (see \cite{dz} 5.2) that for two modular forms $f$ and $g$ of weight $k$ and $l$ the $n$th-Rankin-Cohen Bracket, where $n \geq 0$, given by
\[ (f,g)_n = \sum_{\substack{a,b, \geq 0\\a+b=n}} (-1)^a \binom{k+n-1}{b} \binom{l+n-1}{a} \dif^a f \dif^b g \]
is a modular form of weight $k+l+2n$. In the the case $n>0$ this is a cusp form. For $f= \widetilde{G}_k = \beta_k + [k]$ and $g=\widetilde{G}_l=\beta_l + [l]$ we obtain by using $\dif^a \mb{k}{0} = \frac{(k+a-1)! a!}{(k-1)!} \mb{k+a}{a}$, that
\[ (\widetilde{G}_k , \widetilde{G}_l )_n = \delta_{n,0} \beta_k \beta_l +\gamma^n_{k,l} \cdot C^{2n}_{k,l}\,,\]
with $\gamma^n_{k,l} = \frac{(k-1+n)!}{(k-1)!} \cdot\frac{(l-1+n)!}{(l-1)!} $ and 
\begin{align*}
C^ {2n}_{k,l} &= \beta_k \mb{l+n}{n}+(-1)^n \beta_l \mb{k+n}{n}+ \sum_{\substack{a,b, \geq 0\\a+b=n}} (-1)^a  \mb{k+a}{a} \mb{l+b}{b} \,.
\end{align*}
For all $n\geq1$ and all even $k,l \geq 4$ the function $C^{2n}_{k,l} \in S_k$ is therefore a cusp form of weight $k+l+2n$. This yields a source for relations between bi-brackets since the dimension of $S_k$ is smaller than the possible different $C^{2n}_{k,l}$ .
For example in weight $12$ we have $\dim S_{12} = 1$ and we have the two expressions $\Delta = 12 \cdot 5!^2 \cdot C^4_{4,4} = 5! \cdot 7! \cdot C^2_{4,6}$, with $\Delta = q \prod_{n\geq 0} (1-q^n)^{24}$ being the unique normalized cusp form in this weight. This yields the following relations between bi-brackets
\begin{align*}
7 \mb{5}{1} \cdot \mb{6}{0}-7 \mb{4}{0} \cdot\mb{7}{1}+4 \mb{4}{0}\cdot \mb{6}{2} -2 \mb{5}{1} \cdot \mb{5}{1} = \frac{7}{1440} \mb{7}{1} - \frac{1}{360} \mb{6}{2} +   \frac{1}{8640} \mb{5}{1}  \,.
\end{align*}
%

\section{The spaces $\MD^\ast$ and $\MD^\sh$}

In \cite{H1} it is shown, that every quasi-shuffle Algebra  $(\Q\langle A\rangle,\odot)$  is isomorphic to the shuffle Algebra $(\Q\langle A\rangle,\sh)$. To make this precise define for a composition $i_1+\dots+i_m = n$, where $i_1,\dots,i_m >0$, of a natural number $n$ and a word $w=a_1 a_2  \dots a_n$ the following element in $\Q\langle A\rangle$:
\[ (i_1, \dots, i_m)\{w\} := (a_1 \diamond \dots \diamond a_{i_1}) (a_{i_1+1} \diamond \dots \diamond a_{i_1+i_2}) \dots (a_{i_1+\dots+i_{m-1}+1} \diamond \dots \diamond a_n) \,, \]
where the product is given by the composition of words and $\diamond$ is the product on $\Q A$ belonging to $\odot$. With this define the following two maps
\begin{align*}
 \exp_{\odot}(w) &= \sum_{\substack{1\leq m \leq n\\ i_1+\dots+i_m = n}} \frac{1}{i_1! \dots i_m!}   (i_1, \dots, i_m)\{w\}  \,,\\
\log_{\odot}(w) &=  \sum_{\substack{1\leq m \leq n\\ i_1+\dots+i_m = n}} \frac{(-1)^{n-m}}{i_1 \dots i_m}   (i_1, \dots, i_m)\{w\} \,.
\end{align*}

\begin{prop}(\cite{H1},Thm. 2.5) \label{prop:explog}
 The map $\exp_{\odot}$ is an algebra isomorphism from $(\Q\langle A\rangle,\sh)$ to  $(\Q\langle A\rangle,\odot)$ with the inverse given by $\log_{\odot}$.
\end{prop}
In other words this enables one to give an isomorphism between two arbitrary quasi-shuffle algebras on the same alphabets. We will use this now to define a stuffle version for the brackets and later on the generating series of bi-brackets to define the shuffle brackets. 

Notice that for the brackets, i.e. bi-brackets with $r_1 = \dots = r_l=0$, we also obtain an homomorphism $[.]: (\h^1_z,\boxast) \rightarrow \left( \MD , \cdot \right)$ since we can view $A_z$ as a subset of $A^{\text{bi}}_{z}$. To define the stuffle brackets $[s_1,\dots,s_l]^\ast$, which fulfill the stuffle product, we use the above proposition to deform the quasi-shuffle product $\boxast$ of the brackets into the stuffle product $\ast$, i.e. we use the following compositions of maps to get a algebra homomorphism from $(\h^1_z,\ast)$ to $\MD$. 
\[
\xymatrixrowsep{0.5in}
\xymatrixcolsep{0.7in}
\xymatrix{  
(\h^1_z,\ast)  \ar@{..>}[r]^{[\dots]^\ast} \ar@{->}[d]_{\log_\ast} & (\MD,\cdot)\\ 
 (\h_z^1,\sh)  \ar@{->}[r]_{\exp_\boxast}  & (\h^1_z,\boxast)  \ar@{->}[u]_{[\dots]}
} \]

\begin{dfn}
Define for $s_1,\dots,s_l \in \N$ the \emph{stuffle bracket} $[s_1,\dots,s_l]^\ast$ as the image of $z_{s_1}\dots z_{s_l}$ under the above map, i.e  
\[ [s_1,\dots,s_l]^\ast = [\exp_{\boxast}(\log_{\ast}(z_{s_1}\dots z_{s_l}))] \,. \] 
By $\MD^\ast$ (resp. $\MDA^\ast$) we denote the spaces spanned by all (resp. all with $s_1\geq1$) stuffle brackets and $1$.
\end{dfn}
Remember that the quasi-shuffle product $\boxast$ for brackets was induced by the following map on $\Q A$
\begin{align*}
 z_{s_1} \boxcircle z_{s_2} =&  \sum_{j=1}^{s_1} \lambda^j_{s_1,s_2} z_{j} + \sum_{j=1}^{s_2} \lambda^j_{s_2,s_1}  z_{j}+z_{s_1+s_2} =: z_{s_1+s_2} + \sum_{j\geq 1} \gamma_{s_1,s_2}^j z_j\,,
\end{align*} 
where we define the $\gamma_{s_1,s_2}^j$ just for simplicity of the following formulas. Since $\log_{\ast}(z_{s_1} z_{s_2}) = z_{s_1} z_{s_2} - \frac{1}{2} z_{s_1+s_2}$ and $\exp_{\boxast}(z_{s_1} z_{s_2}) =  z_{s_1} z_{s_2} + \frac{1}{2} z_{s_1+s_2}+\frac{1}{2} \sum_j \gamma_{s_1,s_2}^j z_j $ we obtain $\exp_{\boxast}(\log_{\ast}(z_{s_1} z_{s_2} )) = z_{s_1} z_{s_2} +\frac{1}{2} \sum_j \gamma_{s_1,s_2}^j z_j$, i.e. 
\begin{align*}
[s_1,s_2]^\ast = [s_1,s_2]  + \frac{1}{2} \sum_{j=1}^{s_1} \lambda^j_{s_1,s_2} [j] +  \frac{1}{2} \sum_{j=1}^{s_2} \lambda^j_{s_2,s_1}  [j] \,.
\end{align*}
Similarly one computes the length three case and obtains
\begin{align*}
[s_1,s_2,s_3]^\ast =& [s_1,s_2,s_3] + \frac{1}{2} \sum_{j\geq 0} \gamma_{s_1,s_2}^j [j,s_3] + \frac{1}{2} \sum_{j\geq 0}  \gamma_{s_2,s_3}^j [s_1,j]-\frac{1}{12} \sum_{j\geq 0} \gamma_{s_1+s_2,s_3}^j [j]\\
&- \frac{1}{4} \sum_{j\geq 0} \gamma_{s_1,s_2+s_3}^j [j] + \frac{1}{6} \sum_{j\geq 0} \gamma_{s_1,s_2}^j [j+s_3] +  \frac{1}{6} \sum_{j_1, j_2 \geq 0} \gamma_{s_1,s_2}^{j_1} \gamma_{s_3,j_1}^{j_2} [j_2]\,.
\end{align*}
\begin{ex}For example we have $[1] \cdot [2,1]^\ast = [1,2,1]^\ast + 2[2,1,1]^\ast + [3,1]^\ast + [2,2]^\ast$ with
\begin{align*}
[2,1]^\ast &= [2, 1]-\frac{1}{4} [2],\quad [3,1]^\ast = [3,1] + \frac{1}{24} [2] - \frac{1}{4} [3],\quad  [2,2]^\ast = [2,2]-\frac{1}{12}[2]\,, \\
[2,1,1]^\ast &= [2,1,1]- \frac{3}{4}[2,1]+ \frac{11}{144} [2]-\frac{1}{24} [3]\,, \\
[1,2,1]^\ast &= [1, 2, 1] - \frac{1}{4}[1, 2] - \frac{1}{4} [2, 1] + \frac{1}{72}[2] + \frac{1}{12} [3]\,.
\end{align*} 
\end{ex}
By construction we have the following 
\begin{prop} \label{prop:dimstuff}
Up to lower weight the stuffle brackets equal the brackets and therefore
\[ \dim\left( \grw_k( \MDA)\right)  = \dim\left( \grw_k( \MDA^\ast) \right)\,. \] 
\end{prop}
\begin{prf}
This follow directly from the fact that $\boxast$ and $\ast$ on $\h_z^1$ are equal up to lower weights. 
\end{prf}
In Remark \ref{rem:naivstuffle} we will see that the stuffle brackets can be used to define stuffle regularised the multiple Eisenstein series. However as we will see, even though this version is easy to write down, this will not yield the "best" definition and we will use a more complicated construction. 

We now want to define a $q$-series which is an element in $\bMD$ and which fulfills the "real" shuffle product of multiple zeta values. For $e_1,\dots,e_l \geq 1$ we generalize the generating function of bi-brackets to the following 
\begin{align}\label{def:tribracket} 
\mtt{X_1, & ... & ,\, X_l}{Y_1, & ... & ,\,  Y_l}{e_1, & ... & ,\, e_l} =   \sum_{u_1>\dots>u_l>0} \prod_{j=1}^l E_{u_j}(Y_j) L_{u_j}(X_j)^{e_j} \,. 
\end{align}
So in particular for $e_1 = \dots = e_l=1$ these are the generating functions of the bi-brackets. To show that the coefficients of these series are in $\bMD$ for arbitrary $e_j$ we need to define the differential operator $\mathcal{D}^Y_{e_1,\dots,e_l} := D_{Y_1,e_1}  D_{Y_2,e_2} \dots  D_{Y_l,e_l}$ with
\begin{align*}
D_{Y_j,e} &=  \prod_{k=1}^{e-1} \left(\frac{1}{k}\left(\frac{\partial}{\partial Y_{l-j+1}} - \frac{\partial}{\partial Y_{l-j+2}} \right) - 1\right)\,.
\end{align*}
where we set $\frac{\partial}{\partial Y_{l+1}}=0$. 
\begin{prop}\label{prop:generalpartition}
The coefficients of \eqref{def:tribracket} are in $\bMD$ and it is 
\[\mathcal{D}^Y_{e_1,\dots,e_l}\mt{X_1, \dots, X_l}{Y_1,\dots,Y_l}  = \mtt{Y_1+\dots+Y_l,\, & ... & ,\, Y_1}{X_l,\, X_{l-1}-X_l,\,& ... & ,\,  X_1-X_2}{e_1, & ... & ,\, e_l} \,. \]
\end{prop}
\begin{prf}
By $\frac{\partial}{\partial X} L_n(X) = L_n(X)^2 + L_n(X)$ one inductively obtains 
\[  L_n(Y)^{e+1}  =  \left( \frac{1}{e} \frac{\partial}{\partial Y}  -   1\right) L_n(Y)^e = \prod_{k=1}^{e-1} \left( \frac{1}{k} \frac{\partial}{\partial Y}  -   1 \right) L_n(Y) \,, \]
from which the statement follows after a suitable change of variables. 
\end{prf}
Notice that in the case $e_1=\dots=e_l=1$ this is exactly the partition relation. We now want to define the shuffle brackets $[s_1,\dots,s_l]^\sh$ by using the following well-known fact :
\begin{lem} Let $\mathcal{A}$ be an algebra spanned by elements $a_{s_1,\dots,s_l}$ with $s_1,\dots, s_l \in \N$, let $H(X_1,\dots,X_l) = \sum_{s_j} a_{s_1,\dots,s_l} X_1^{s_1-1} \dots X_1^{s_l-1}$ be the generating functions of these elements and define for $f \in \Q[[X_1,\dots,X_l]]$
\[ f^{\sharp}(X_1,\dots,X_l) = f(X_1+\dots+X_l, X_2+\dots+X_l, \dots, X_l) \,. \] 
Then the following two statements are equivalent 
\begin{enumerate}[i)]
\item The map $(\h^1_{xy}, \sh) \rightarrow \mathcal{A}$ given by $z_{s_1} \dots z_{s_j} \mapsto a_{s_1,\dots,s_l}$ is an algebra homomorphism. 
\item For all $r,s \in \N$  it is
\[ H^{\sharp}(X_1,\dots,X_r) \cdot H^{\sharp}(X_{r+1}, \dots, X_{r+s}) = H^{\sharp}(X_1,\dots,X_{r+s})_{\vert sh_r^{(r+s)}}\,, \]
where $sh_r^{(r+s)} = \sum_{\sigma\in\Sigma(r,s)}\sigma$ in the group ring $\Z[\mathfrak{S}_{r+s}]$ and the symmetric group $\mathfrak{S}_r$ acts on $\Q[[X_1,\ldots,X_r]]$ by $(f\big|\sigma)(X_1,\ldots,X_r)= f(X_{\sigma^{-1}(1)},\ldots,X_{\sigma^{-1}(r)})$\,.
\end{enumerate}
\end{lem}
\begin{prf} A Lemma of this type was used in \cite{IKZ}. \end{prf}

\begin{thm} \label{thm:shufflebracket} For $s_1,\dots,s_l\in \N$ define $[s_1,\dots,s_l]^\sh \in \bMD$ as the coefficients of the following generating function
\begin{align*}
&H_{\sh}(X_1,\dots,X_l)  = \sum_{s_1,\dots,s_l \geq 1} [s_1,\dots,s_l]^\sh X_1^{s_1-1} \dots X_l^{s_l-1} \\
&:=\sum_{\substack{1\leq m \leq l\\i_1 + \dots + i_m = l}}\frac{1}{i_1! \dots i_m!} \mathcal{D}^Y_{i_1,\dots,i_m} \mt{X_1,X_{i_m+1},X_{i_{m-1}+i_m+1}, \dots, X_{i_2+\dots+i_m+1}}{Y_1,\dots,Y_l} _{\big\vert Y=0} \,.
\end{align*}
Then we have the following to statements
\begin{enumerate}[i)]
\item The $[s_1,\dots,s_l]^\sh$ fulfill the shuffle product, i.e. 
\[ H^{\sharp}_{\sh}(X_1,\dots,X_r) \cdot H^{\sharp}_{\sh}(X_{r+1}, \dots, X_{r+s}) = H^{\sharp}_{\sh}(X_1,\dots,X_{r+s})_{\vert sh_r^{(r+s)}} \,. \]
\item For $s_1\geq 1,\, s_2,\dots,s_l \geq 2$ we have $[s_1,\dots,s_l]^\sh = [s_1,\dots,s_l]$.
\end{enumerate}
\end{thm}
\begin{prf} The first part of the proof is basically the same as in the discussion in section 4.1 in \cite{BT} but with a reverse order and some changes in the notation.
Consider the alphabet $A=\left\{ \tbinom{y}{n} \mid n \in \N \,, y \in Y_\Z  \right\}$, where $Y_{\Z}$ is the set of finite sums of the elements in $Y=\{Y_1, Y_2,\dots\}$. We denote a word in these letters by $\tbinom{y_1,\dots,y_l}{n_1,\dots,n_l}$. For two letters $a, b \in A$ define $a \diamond b \in A$ as the component-wise sum. With this we can equip $\Q\langle A \rangle$ with the quasi-shuffle product $\odot$ \eqref{qshp}  and therefore obtain a quasi-shuffle algebra $( \Q\langle A \rangle , \odot )$. It is easy to see that the map $( \Q\langle A \rangle , \odot ) \rightarrow \bMDG$ given by
\[ \tbinom{y_1,\dots,y_l}{n_1,\dots,n_l} \longmapsto  \mtt{0, & ... & ,\, 0}{y_1, & ... & ,\,  y_l}{n_1, & ... & ,\, n_r} \] 
is an algebra homomorphism. Using now Proposition \ref{prop:explog} the series $h$ defined by the exponential map
\[  h(X_1,\dots,X_r) =\sum_{\substack{1\leq m \leq n\\ i_1+\dots+i_m = n}} \frac{1}{i_1! \dots i_m!} \mtt{0, & ... & ,\, 0}{Y_1, & ... & ,\,  Y_m}{i_1, & ... & ,\, i_m} \,, \] 
where  $Y_j = X_{i_1+\dots+i_{j-1}} + \dots + X_{i_1 + \dots + i_j}$ with $X_0:=0$, fulfills the (index-)shuffle product i.e. 
\[ h(X_1,\dots,X_r) \cdot h(X_{r+1}, \dots, X_{r+s}) = h(X_1,\dots,X_{r+s})_{\vert sh_r^{(r+s)}} \,.\]
We now set $H_{\sh}(X_1,\dots,X_l) := h(X_l,X_{l-1}-X_l,\dots, X_1-X_2)$ and by the same argument as in Theorem 4.3 in \cite{BT} it is 
\[ H^{\sharp}_{\sh}(X_1,\dots,X_r) \cdot H^{\sharp}_{\sh}(X_{r+1}, \dots, X_{r+s}) = H^{\sharp}_{\sh}(X_1,\dots,X_{r+s})_{\vert sh_r^{(r+s)}} \,. \]
Combining the definition of $h$ and $H_{\sh}$ we observe that $H_{\sh}(X_1,\dots,X_r)$ equals
\[ \sum_{\substack{1\leq m \leq n\\ i_1+\dots+i_m = n}} \frac{1}{i_1! \dots i_m!} \mtt{0, & ... & ,\, 0}{X_{r-i_1+1},\, X_{r-i_1-i_2+1} - X_{r - i_1 + 1} , &...& ,\,  X_1 - X_{r-i_1-\dots-i_{m-1}+1}}{i_1, & ... & ,\, i_m} \,.\] 
We now apply Proposition \ref{prop:generalpartition} to this and obtain i) of the Theorem. To prove ii) one checks that the only summand on the right hand side, where {\bf all} variables $X_2,\dots,X_l$ appear, is the one with $i_1= \dots = i_m=1$ which is exactly $[s_1,\dots,s_l] X^{s_1-1} \dots X_l^{s_l-1}$. Therefore the shuffle bracket $[s_1,\dots,s_l]^\sh$ where $s_2,\dots,s_l\geq 2$ is given by the bracket $[s_1,\dots,s_l]$. 
\end{prf}

For low length we obtain the following examples: 
\begin{cor}\label{cor:explicitshufflebracket}
It is $[s_1]^\sh = [s_1]$ and for $l=2,3,4$  the $[s_1,\dots,s_l]^\sh$are given by
\begin{enumerate}[i)]
\item $\begin{aligned}[t]
[s_1,s_2]^\sh &= [s_1,s_2] + \delta_{s_2,1}\cdot \frac{1}{2} \left( \mb{s_1}{1} - [s_1] \right) \,,\\
\end{aligned}$
\item $\begin{aligned}[t]
[s_1,s_2,s_3]^\sh &= [s_1,s_2,s_3]+ \delta_{s_3,1}\cdot \frac{1}{2} \left( \mb{s_1,s_2}{0,1} - [s_1,s_2] \right)\\
&+\delta_{s_2,1}\cdot \frac{1}{2}\left(  \mb{s_1,s_3}{1,0} - \mb{s_1,s_3}{0,1}- [s_1,s_3] \right)\\
&+\delta_{s_2 \cdot s_3,1}\cdot \frac{1}{6}\left(  \mb{s_1}{2} - \frac{3}{2}\mb{s_1}{1}+ [s_1] \right) \,,\\
\end{aligned}$
\item $\begin{aligned}[t]
[s_1,s_2,s_3&,s_4]^\sh = [s_1,s_2,s_3,s_4] + \delta_{s_4,1}\cdot \frac{1}{2} \left( \mb{s_1,s_2,s_3}{0,0,1} - [s_1,s_2,s_3] \right) \\
+ \delta_{s_3,1} &\cdot\frac{1}{2} \left( \mb{s_1,s_2,s_4}{0,1,0}-\mb{s_1,s_2,s_4}{0,0,1} + [s_1,s_2,s_4] \right)\\
+ \delta_{s_2,1} &\cdot\frac{1}{2} \left( \mb{s_1,s_3,s_4}{1,0,0}-\mb{s_1,s_3,s_4}{0,1,0} + [s_1,s_3,s_4] \right)\\
+ \delta_{s_2 \cdot s_4,1} &\cdot\frac{1}{4} \left( \mb{s_1,s_3}{1,1}-2\mb{s_1,s_3}{0,2}-\mb{s_1,s_3}{1,0} + [s_1,s_3] \right)\\
+ \delta_{s_3 \cdot s_4,1} &\cdot\frac{1}{6} \left( \mb{s_1,s_2}{0,2}-\frac{3}{2}\mb{s_1,s_2}{0,1} + [s_1,s_2] \right)\\
+ \delta_{s_2 \cdot s_3,1} &\cdot\frac{1}{6} \left( \mb{s_1,s_4}{0,2}-\mb{s_1,s_4}{1,1}+\frac{3}{2}\mb{s_1,s_4}{0,1}+\mb{s_1,s_4}{2,0}-\frac{3}{2}\mb{s_1,s_4}{1,0} + [s_1,s_4] \right)\\
+ \delta_{s_2 \cdot s_3 \cdot s_4,1} &\cdot\frac{1}{24} \left( \mb{s_1}{3}-2\mb{s_1}{2}+\frac{11}{6} \mb{s_1}{1} - [s_1] \right) \,.
\end{aligned}$
\end{enumerate}
\end{cor}
\begin{prf}
This follows by calculating the coefficients of the series $G_\sh$ in Theorem \ref{thm:shufflebracket}.
\end{prf}

\begin{prop}\label{prop:len2conj}
For all $s_1,s_2 \geq 1$ it is
\begin{align*}
\mb{s_1,s_2}{1,0}, \mb{s_1,s_2}{0,1} \in \operatorname{Fil}^{\operatorname{W},\operatorname{L}}_{s_1+s_2+1,3}(\MD)
\end{align*}
\end{prop}
\begin{prf}
 First notice that from $\mb{s_1,s_2}{1,0} \in \MD$ by the stuffle product for bi-brackets $\mb{s_1}{1} \cdot [s_2]$ one deduces $\mb{s_2,s_1}{0,1} \in \MD$. Since the shuffle brackets fulfill the shuffle product we have
\begin{align*}
[s_1,s_2]^\sh \cdot [1] = 2 [s_1,s_2,1]^\sh + 2[s_1,1,s_2]^\sh + 2 [1,s_1,s_2]^\sh + \sum_{a,b,c \geq 2} \nu_{a,b,c} [a,b,c]^\sh
\end{align*}
for some $\nu_{a,b,c}  \in \Q$. By Proposition \ref{cor:explicitshufflebracket} the brackets $[s_1,s_2]^\sh,\, [1,s_1,s_2]^\sh$ and $[a,b,c]^\sh$ with $a,b,c\geq2$ are elements of $\MD$, i.e.  $2 [s_1,s_2,1]^\sh + 2[s_1,1,s_2]^\sh \in \MD$. Using the explicit formula for the length three shuffle brackets it is easy so check that 
\[ 2 [s_1,s_2,1]^\sh + 2[s_1,1,s_2]^\sh = \left\{
\begin{array}{cl} \mb{s_1,s_2}{1,0} \,, & s_2 > 1,  \\ 2 \mb{s_1,1}{0,1} \,, & s_2=1
\,. \end{array}  \right. \mod \MD \,,\]
which proves the statement. 
\end{prf}
%

Finally we give some numerical results on the dimension of the space spanned by the shuffle brackets $[s_1,\dots,s_l]^\sh$. Denote by $\MD^\sh$ the $\Q$-vector space spanned by all $[s_1,\dots,s_l]^\sh$ and $1$ and $\MDA^\sh$ spanned by those where $s_1>1$. By the use of the computer the author was able to give lower bounds for the dimension of $\grw_k( \MD^\sh)$ for $k\leq 10$ by using a fast implementation of the bi-brackets in Pari GP
\begin{table}[H]\footnotesize
\begin{center}
\begin{tabular} { c| c| c | c | c |c|c|c|c|c|c|c|c|c|c|c|}
$k$ 		&	0 & 1	&	2 & 3 & 4	& 5 & 6 &  7 & 8 & 9 & 10  \\ \hline
$\dim\left(\grw_k( \MDA^\sh) \right) \geq$  & 1 & 0 & 1 & 2 & 3 & 6 & 10 & 18 & 32 & 56 & 100 	\\ 
\end{tabular} 
\caption{Lower bounds for $\dim\left(\grw_k( \MDA^\sh) \right)$.} 
\end{center}
\end{table}
We observe that these numbers coincide with the conjectured dimension for $\grw_k( \MDA)$ given in \cite{BK}. Setting $d'_0=1$, $d'_1=0$, $d'_2=1$,  $d'_3=2$, $d'_4=3$
and for $k\geq5$:
 \begin{align*} 
d'_k = 2 d'_{k-2} + 2 d'_{k-3} \,,
\end{align*}
we have the following conjecture.
\begin{conj}
The dimensions of  $\grw_k( \MDA)$ and  $\grw_k( \MDA^\sh)$ coincide and they are given by
\[ \dim\left( \grw_k( \MDA)\right) = \dim\left( \grw_k( \MDA^\sh) \right) = d'_k \,. \]
\end{conj}
Recall that by Proposition \ref{prop:dimstuff} this conjecture would also imply  $\dim\left( \grw_k( \MDA^\ast) \right)  = \dim\left( \grw_k( \MDA^\sh) \right)$. 
\begin{rem}
In the case of multiple zeta values the shuffle product is an easy consequence of the expression as an iterated integral. It is therefore a natural question whether there is also some kind of iterated integral expression from which the shuffle product follows. This was done for other $q$-analogue models of MZV in \cite{Zh} and \cite{JMK} by the use of iterated Jackson integrals. 
\end{rem}

\section{Multiple Eisenstein series $G$, $G^\sh$ and $G^\ast$}
In \cite{BT} the authors defined regularized multiple Eisenstein series via the use of the coproduct structure on the space of formal iterated integrals. We will recall the basic facts in the following. Since in \cite{BT} a {\bf different order} in the definition of MZV was used we will use the following definitions of MZV and MES for the rest of this section:  
\begin{align*}
 \zetarev(s_1,\ldots,s_l)&=\sum_{{\bf 0<n_1<\cdots<n_r}} \frac{1}{n_1^{s_1}\cdots n_r^{s_l}} \quad (s_1,\ldots,s_{l-1} \geq 1, s_l \geq 2) \,,
\end{align*}
obviously we have $ \zetarev(s_1,\ldots,s_l) =  \zeta(s_l,\ldots,s_1)$.
For the multiple Eisenstein series we will use a similar notation, i.e. the multiple Eisenstein series in \cite{BT} are denoted $\grev$ here and the relation to the multiple Eisenstein series $G$ given in the introduction is 
\[\grev_{s_1,\dots,s_l}(\tau) = G_{s_l,\dots,s_1}(\tau) \] 
and similarly $\gshrev$ and $\gstrev$. The notation for multiple Eisenstein series $G$ coincides with the original paper \cite{gkz} and the work \cite{Bach} on multiple Eisenstein series but differs from \cite{BT}.

\begin{dfn}
For integers $s_1,\ldots,s_{l-1}\ge2$ and $s_l \geq 3$, we define the \emph{multiple Eisenstein series} $\grev_{s_1,\ldots,s_l}(\tau)$ on $\mathbb{H}$ by 
\[ \grev_{s_1,\ldots,s_l}(\tau) =   \sum_{0\prec \lambda_1\prec \cdots\prec \lambda_r} \frac{1}{\lambda_1^{s_1}\cdots \lambda_r^{s_l}} \,, \]
where $\lambda_i \in \Z \tau+ \Z$ are lattices points and the order $\prec$ on $\Z + \Z \tau$ is given by 
\[ m_1 \tau + n_1 \prec m_2 \tau + n_2 :\Leftrightarrow \left( m_1 < m_2 \vee (m_1 = m_2 \wedge n_1 < n_2) \right)\,.\]
\end{dfn}
\begin{rem}\label{rem:stuffleprod}
It is easy to see that these are holomorphic functions in the upper half plane and that they fulfill the stuffle product, i.e. it is for example
\[ \grev_3(\tau) \cdot \grev_4(\tau) = \grev_{4,3}(\tau) + \grev_{3,4}(\tau) + \grev_{7}(\tau)\,. \] 
The condition $s_l \geq 3$ is necessary for absolutely convergence of the sum. By choosing a specific way of summation we can also restrict this condition to get a definition of $\grev_{s_1,\ldots,s_l}(\tau)$ with $s_l=2$ which also satisfies the stuffle product (see \cite{BT} for detail). 
\end{rem}

Recall that we denote by $\MZB \subset \C[[q]]$ the space spanned by all $q$-series given by products of MZV, powers of $(-2 \pi i)$ and bi-brackets. In \cite{Bach} the Fourier expansion of multiple Eisenstein series was calculated. In particular the results in \cite{Bach} show that we can consider $G_{s_1,\ldots,s_l}$ as well as $\grev_{s_1,\ldots,s_l}$ to be elements in $\MZB$ by setting $q=e^{2\pi i\tau}$. For example 
\begin{align*}
\grev_{2,3}(\tau) &=\zetarev(2,3) + 3 \zetarev(3) g_2(q)  + 2 \zetarev(2) g_3(q) + g_{2,3}(q) \in \MZB \,,
\end{align*}
where for all $s_1,\dots,s_l \geq 1$ we write $g_{s_1,\ldots,s_l}(q) = (-2\pi i)^{s_1+\dots+s_l}[s_l,\ldots,s_1]$. We will also use the following notation 
\begin{align*}
g^\sh_{s_1,\ldots,s_l}(q) &= (-2\pi i)^{s_1+\dots+s_l}[s_l,\ldots,s_1]^\sh \,, \\
 g\tbinom{s_1,\dots,s_l}{r_1,\dots,r_l}(q) &= (-2\pi i)^{s_1+r_1+\dots+s_l+r_l} \mb{s_l,\dots,s_1}{r_l,\dots,r_1} \,.
\end{align*}
Later we will suppress the dependence of $q$ and $\tau$ and just write $g_{s_1,\ldots,s_l}$ instead of $g_{s_1,\ldots,s_l}(q)$ and similar for the other functions considered above.

Following Goncharov (\cite{G}) the authors in \cite{BT} consider the algebra ${\mathcal I}$ generated by the elements $\mathbb{I}(a_0;a_1,\ldots,a_N;a_{N+1})$, where $a_i\in\{0,1\}, N\ge0$, with the product given by the shuffle product $\sh$ together with relations coming from real iterated integrals (see \cite{G} and \cite{BT} for details). This space has the structure of a Hopf algebra with the coproduct given by
\begin{align} \label{deltag}
\begin{split}
 &\Delta_G \left( \mathbb{I}(a_0;a_1,\ldots,a_N;a_{N+1}) \right) :=\\
 & \sum \Big(\prod_{p=0}^k \mathbb{I}(a_{i_p};a_{i_p+1},\ldots,a_{i_{p+1}-1};a_{i_{p+1}}) \Big) \otimes  \mathbb{I}(a_0;a_{i_1},\ldots,a_{i_k};a_{N+1}),
 \end{split}
\end{align}
where the sum runs over all $i_0=0<i_1<\cdots<i_k<i_{k+1}=N+1$ with $0\le k \le N$. The triple $({\mathcal I},\sh, \Delta_G)$ is a commutative graded Hopf algebra over $\Q$. For integers $n\ge0,s_1,\ldots,s_l\ge1$, we set 
\[I_{n}(s_1,\ldots,s_l):=I(0;\underbrace{0,\ldots,0}_{n},\underbrace{1,0,\ldots,0}_{s_1},\ldots,\underbrace{1,0,\ldots,0}_{s_l};1).\]
In particular, we write $I(s_1,\ldots,s_l)$ to denote $I_0(s_1,\ldots,s_l)$. The quotient space $\mathcal{I}^1=\mathcal{I}/\mathbb{I}(0;0;1)\mathcal{I}$ also has the structure of a Hopf algebra with the same coproduct and due to Proposition 3.2 in \cite{BT} the elements $I(s_1,\ldots,s_l)$ form a basis of $\mathcal{I}^1$, i.e. as a $\Q$-algebra the space $\mathcal{I}^1$ is isomorphic to $(\h_{xy}^1,\sh)$ by sending $I(s_1,\ldots,s_l)$ to $z_{s_1} \dots z_{s_l}$ (Notice that since the order changed, we write $z_j = y x^{j-1}$ in this section). In the following we therefore consider $\h^1_{xy}$ as a Hopf algebra with the above coproduct. 

\begin{prop}{\cite{IKZ}}(shuffle \& stuffle regularised MZV)
There exist algebra homomorphisms $Z^{\sh} : (\h_{xy}^1 , \sh) \rightarrow \MZ$ and $Z^{\ast} : (\h^1_z, \ast) \rightarrow \MZ$ with $\zetashrev(s_1,\dots,s_l) = Z^\sh(z_{s_1} \dots z_{s_l})$ and $\zetashrev(s_1,\dots,s_l) = Z^\sh(z_{s_1} \dots z_{s_l})$ such that
\[ \zetashrev(s_1,\ldots,s_l) = \zetastrev(s_1,\ldots,s_l) = \zetarev(s_1,\ldots,s_l) \] 
for $s_1,\dots,s_{l-1} \geq 1$ and $s_l \geq 2$. They are uniquely determined by $Z^{\sh}(z_1)=Z^{\ast}(z_1)=0$. 
\end{prop}
\begin{prf}
This follows from the results of section $2$ in \cite{IKZ}. 
\end{prf}

 We now recall the definition of $\gshrev$ from \cite{BT}. 
\begin{dfn}
For integers $s_1,\ldots,s_l\ge1$, define the $q$-series $\gshrev_{s_1,\ldots,s_l}(q)\in \MZB$, called  \emph{(shuffle) regularized multiple Eisenstein series}, as 
\[ \gshrev_{s_1,\ldots,s_l}(q) := m\left( (Z^{\sh}\otimes \mathfrak{g}^{\sh})\circ \Delta_G \big(z_{s_1} \dots z_{s_l} \big)\right)\,,\]
where  $\mathfrak{g}^{\sh} : (\h_{xy}^1, \sh) \rightarrow \C[[q]]$ is the algebra homomorphism defined by $\mathfrak{g}^{\sh}(z_{s_1} \dots z_{s_l}) = g^\sh_{s_1,\ldots,s_l}(q)$ and $m$ denotes the multiplication given by $m: a \otimes b \mapsto a\cdot b$.
\end{dfn}
We can view $\gshrev $ as an algebra homomorphism $\gshrev : (\h_{xy}^1, \sh) \rightarrow \MZB$ such that the following diagram commutes
\[
\xymatrix{  
(\h_{xy}^1,\sh)  \ar@{->}[r]^-{\Delta_G} \ar@{->}[d]_{\gshrev} & (\h_{xy}^1,\sh)  \otimes (\h_{xy}^1,\sh)  \ar@{->}[d]^{Z^{\sh} \otimes\, \mathfrak{g}^{\sh} }  \\ 
 \MZB &  \MZ \otimes \C[[q]]  \ar@{->}[l]^-{m}\\
} \]

Summarizing the results of \cite{BT} we have
\begin{thm}{\cite{BT}}\label{thm:restrds}
For all $s_1,\ldots,s_l\ge1$ and $q=e^{2\pi i\tau}$ with $\tau \in \Ha$ the regularised multiple Eisenstein series $\gshrev_{s_1,\ldots,s_l}(q)$ have the following properties:
\begin{enumerate}[i)]
\item They are holomorphic functions on the upper half plane having a Fourier expansion with the regularised multiple zeta values as the constant term. 
\item They fulfill the shuffle product, i.e. we have an algebra homomorphism $(\h_{xy}^1 , \sh) \rightarrow \MZB$ by sending the generators $z_{s_1} \dots z_{s_l}$ to $\gshrev_{s_1,\ldots,s_l}(q)$.
\item For integers $s_1,\ldots,s_l\ge2$ they equal the multiple Eisenstein series
\[ \gshrev_{s_1,\ldots,s_l}(q)=\grev_{s_1,\ldots,s_l}(q) \]
and therefore they fulfill the stuffle product (see Remark \ref{rem:stuffleprod}) in these cases. 
\end{enumerate}
\end{thm}
Theorem \ref{thm:restrds} provides a large family of linear relations between the $\gshrev$, since one can write the product $\gshrev_{s_1,\dots,s_l}(q) \cdot \gshrev_{r_1,\dots,r_m}(q)$ in two different ways whenever one has $s_1,\dots,s_l,r_1,\dots,r_m \geq 2$ by using the stuffle and shuffle product formula. We will call these relations the restricted double shuffle relations, since they are just a subset of all (finite) double shuffle relations of MZV, where the indices $s_j$ and $r_i$ are additionally allowed to be $1$ whenever $j<l$ and $i<m$. We compare the number of both relations at the end of this paper. 

Numerical experiments suggest (see the dimension discussion at the end of \cite{BT}), that there are additional relations between the $\gshrev$ coming from the double shuffle relations where some indices are also allowed to be $1$. It is therefore interesting to understand the exact failure of the stuffle product for the regularised multiple Eisenstein $\gshrev$ which seems not to be covered best possible by the Theorem \ref{thm:restrds}. In the following we want to sketch a possible approach to answer this question. The basic idea is to define stuffle regularised  multiple Eisenstein series $\gstrev_{s_1,\dots,s_l}$  which equals the shuffle regularised ones in most of the cases. For this we need the following: For an arbitrary quasi-shuffle algebra $\Q\langle A \rangle$ define on  the following coproduct for a word $w$ 
\[ \Delta_H(w) = \sum_{u v = w} u \otimes v\,. \] 
Then it is known due to Hoffman (\cite{H1}) that the space $\left( \Q\langle A \rangle, \odot , \Delta_H \right)$ has the structure of a bialgebra. With this we try to mimic the definition of the $\gshrev$ and use the coproduct structure on the space $(\h^1_z, \ast, \Delta_H)$ to define $\gstrev$, i.e. we consider the following diagram
\[
\xymatrix{  
 (\h^1_z,\ast)  \ar@{->}^-{\Delta_H} [r]\ar@{->}[d]_{\gstrev}  &(\h^1_z,\ast) \otimes (\h^1_z,\ast) \ar@{->}[d]^{Z^{\ast} \otimes\, \mathfrak{g}^{\ast} }   \\ 
 \C[[q]] & \MZ \otimes \C[[q]] \ar@{->}[l]^-{m}   &\\
} \]
with a suitable choice of an algebra homomorphism $\mathfrak{g}^{\ast} : (\h^1_z, \ast) \rightarrow \C[[q]]$.

\begin{rem}\label{rem:naivstuffle}
One naive way to define $\mathfrak{g}^\ast$ would be to define it on the generator $w=z_{s_l} \dots z_{s_l}$ by $(-2\pi i)^{s_1+\dots+s_l}[s_l,\ldots,s_1]^\ast\,$ which would yield stuffle regularised the multiple Eisenstein series which coincide with the $\gshrev$ in the length one case. But already in length two this differs from the original multiple Eisenstein series even when all $s_j \geq 2$ for example it is
\begin{align*}
\grev_{2,3}(\tau) = \grev_{2,3}^\sh(\tau) &=\zetarev(2,3) + 3 \zetarev(3) g_2(q)  + 2 \zetarev(2) g_3(q) + g_{2,3}(q)
\end{align*}
but the naive approach would give $\zetarev(2,3) + 2 \zetarev(2) g_3(q) + g_{2,3}(q)$. Even though these are similar this seems not to be the definition we want and we need to find an alternative definition for $\mathfrak{g}^\ast$ in the following such that $\gstrev$ coincide with the original multiple Eisenstein series.
\end{rem}

Motivated by the calculation of the Fourier expansion of multiple Eisenstein series described in $\cite{Bach}$ and $\cite{BT}$ we consider the following construction. 
%
\begin{constr} \label{const}
Given a $\Q$-algebra $(A,\cdot)$ and a family of homomorphism
\[ \left\{ w \mapsto f_w(m) \right\}_{m \in \N} \]
 from  $(\h^1_z, \ast)$ to $(A,\cdot)$, we define for $w \in \h^1_z$ and $M \in \N$ 
\[ F_w(M) := \sum_{\substack{1 \leq k \leq l(w)\\w_1 \dots w_k = w\\0< m_1 < \dots <m_k <M}}  f_{w_1}(m_1) \dots f_{w_k}(m_k) \in A \,, \]
where $l(w)$ denotes the length of the word $w$ and $w_1 \dots w_k = w$ is a decomposition of $w$ into $k$ words in $\h^1_z$. 
\end{constr}
\begin{prop}\label{prop:construction}
For all $M \in \N$ the map from $(\h^1_z, \ast)$ to $(A,\cdot)$ defined by $w \mapsto F_w(M)$ is an algebra homomorphism, i.e. $\left\{ w \mapsto F_w(m) \right\}_{m \in \N}$ is again a family  of homomorphism as in the Construction \ref{const}. 
\end{prop}
\begin{prf}
We use the coproduct structure on $\left( \h^1_z , \ast , \Delta_H \right)$ to prove the statement by induction over $M$. It is $F_w(1) =0$ which clearly fulfills the stuffle product. For the induction step one checks that $F_w(M+1) = \sum_{\substack{uv =w}}  F_u(M) f_v(M)$ which is exactly the image of $w$ under $(F(M) \otimes  f(M)) \circ \Delta_H$, i.e. it fulfills the stuffle product by the induction hypothesis. 
\end{prf}

For a word $w=z_{s_1} \dots z_{s_l} \in \h^1$ we also write in the following $f_{s_1,\dots,s_l}(m):=f_w(m)$ and similarly $ F_{s_1,\dots,s_l}(M):=F_w(M)$. 

\begin{ex} Let $f_w(m)$ be as in the construction. In small lengths the $F_w$ are given by
\[ F_{s_1}(M) = \sum_{0< m_1 <M} f_{s_1}(m_1) \,, \quad F_{s_1,s_2}(M) = \sum_{0< m_1 <M} f_{s_1,s_2}(m_1) + \sum_{0< m_1 < m_2<M} f_{s_1}(m_1)  f_{s_2}(m_2) \,\]
and one can check directly by the use of the stuffle product for the $f_w$ that 
\begin{align*}
&F_{s_1}(M) \cdot F_{s_2}(M) =  \sum_{0< m_1 <M} f_{s_1}(m_1) \cdot  \sum_{0< m_2 <M} f_{s_2}(m_2) \\
&=   \sum_{0< m_1 < m_2<M} f_{s_1}(m_1)  f_{s_2}(m_2) +  \sum_{0< m_2 < m_1<M}   f_{s_2}(m_2) f_{s_1}(m_1) +  \sum_{0<  m_1> M}   f_{s_1}(m_1) f_{s_2}(m_1)  \\
&=  \sum_{0< m_1 < m_2<M} f_{s_1}(m_1)  f_{s_2}(m_2) +  \sum_{0< m_2 < m_1<M}   f_{s_2}(m_2) f_{s_1}(m_1) \\
&+ \sum_{0< m_1 <M} \left(f_{s_1,s_2}(m_1) +  f_{s_2,s_1}(m_1) + f_{s_1+s_2}(m_1) \right)  \\
&=  F_{s_1,s_2}(M)+ F_{s_2,s_1}(M) + F_{s_1+s_2}(M)  \,.
\end{align*}
\end{ex}

Let us now give an explicit example for maps $f_w$ in which we are interested. For this we need to define the following
\begin{dfn}\label{def:multitangent}
For integers $s_1,\dots,s_l \ge 1$ with $s_1,s_l\ge2$ we define a holomorphic function $\Psi_{s_1,\ldots,s_l}(z)$ on $\C- \Z$ called the \emph{multitangent function} by
\[ \Psi_{s_1,\ldots,s_l}(z) = \sum_{\substack{n_1<\cdots<n_l\\ n_j \in \Z}} \frac{1}{(z+n_1)^{s_1}\cdots (z+n_l)^{s_l}}.\]
When $l=1$ we refer to $\Psi_{s_1}(z)$ as the \emph{monotangent function}.
\end{dfn}

In \cite{Boulliot} the author uses the notation $\mathcal{T}e^{n_r,\ldots,n_1}(z)$ which corresponds to our $\Psi_{n_1,\ldots,n_r}(z)$ and showed that the series defining $\Psi_{n_1,\ldots,n_r}(z)$ converges absolutely when $n_1,\ldots,n_r\ge2$. These functions fulfill (for the cases they are defined) the stuffle product. The multitangent functions appear in the calculation of the Fourier expansion of the multiple Eisenstein series $\grev_{s_1,\dots,s_l}$ (see \cite{Bach}, \cite{BT}), for example in length two it is
\[ \grev_{s_1,s_2}(\tau) = \zetarev(s_1,s_2) + \zetarev(s_1) \sum_{m_1 > 0} \Psi_{s_2}(m_1 \tau) + \sum_{m_1>0} \Psi_{s_1,s_2}(m_1 \tau) +  \sum_{m_1 > m_2 > 0} \Psi_{s_1}(m_1\tau) \Psi_{s_2}(m_2\tau) \,.\] 
One nice result of \cite{Boulliot} is a regularization of the multitangent function to get a definition of $\Psi_{s_1,\ldots,s_l}(z)$ for all $s_1,\dots,s_l \in \N$.  We will use this result together with the above construction to recover the Fourier expansion of the multiple Eisenstein series. 
\begin{thm}(\cite{Boulliot})\label{thm:olivier}
For all $s_1,\dots,s_l \in \N$  there exist holomorphic functions $\Psi_{s_1,\dots,s_l}$ on $\Ha$ with the following properties
\begin{enumerate}[i)]
\item Setting $q=e^{2\pi i \tau}$ for $\tau \in \Ha$ the map  $w \mapsto \Psi_w(\tau)$ defines an algebra homomorphism from $(\h^1_z, \ast)$ to $(\C[[q]],\cdot)$.
\item In the case $s_1,s_l \geq 2$ the  $\Psi_{s_1,\dots,s_l}$ are given by the multitangent functions in Definition \ref{def:multitangent}.
\item The monotangents functions have the $q$-expansion given by
\[ \Psi_1(\tau) = \frac{\pi}{\tan(\pi \tau)} =  (-2 \pi i)\left(\frac{1}{2} + \sum_{n>0} q^n \right),\quad  \Psi_k(\tau) =  \frac{(-2\pi i)^k}{(k-1)!} \sum_{n>0} n^{k-1} q^n \, \text{ for } k\geq 2.\] 
\item (Reduction into monotangent function) Every $\Psi_{s_1,\dots,s_l}(\tau)$ can be written as a $\MZ$-linear combination of monotangent functions. There are explicit $\epsilon^{s_1,\dots,s_l}_{i,k} \in \MZ$ s.th. 
\[\Psi_{s_1,\dots,s_l}(\tau) = \delta^{s_1,\dots,s_l} + \sum_{i=1}^l \sum_{k=1}^{s_i} \epsilon^{s_1,\dots,s_l}_{i,k} \Psi_k(\tau) \,,\]
where $ \delta^{s_1,\dots,s_l} = \frac{(\pi i)^l}{l!}$ if $s_1=\dots=s_l=1$ and $l$ even and $\delta^{s_1,\dots,s_l} = 0$ otherwise. 
For $s_1>1$ and $s_l>1$ the sum on the right starts at $k=2$, i.e. there are no $\Psi_1(\tau)$ appearing and therefore there is no constant term in the $q$-expansion.
\end{enumerate}
\end{thm}
\begin{prf}
This is just a summary of the results in Section $6$ and $7$ of \cite{Boulliot}. The last statement is given by Theorem 6 there. 
\end{prf} 
Due to iv) in the Theorem the calculation of the Fourier expansion of multiple Eisenstein series, where ordered sums of multitangent functions appear, reduces to ordered sums of monotangent functions. The connection of these sums to the brackets, i.e. to the functions $g$, is given by the following fact which can be seen by using iii) of the above Theorem. For $n_1,\dots,n_r \geq 2$ it is
\begin{align*}
&g_{s_1,\ldots,s_r}(q) = \sum_{0 < m_1 <\dots <m_l} \Psi_{s_1}(m_1\tau) \dots \Psi_{s_l}(s_l\tau) \,.
\end{align*}
For $w \in \h^1_z$ we now use the Construction \ref{const} with $A=\C[[q]]$ and the family of homomorphism $\{ w \mapsto \Psi_w(n \tau) \}_{n\in \N}$ to define
\[ \mathfrak{g}^{\ast,M}(w) := (-2\pi i)^{|w|} \sum_{\substack{1 \leq k \leq l(w)\\w_1 \dots w_k = w}} \sum_{0< m_1 < \dots <m_k <M} \Psi_{w_1}(m_1 \tau) \dots \Psi_{w_k}(m_k \tau)  \,. \] 
From Proposition \ref{prop:construction} and the Theorem \ref{thm:olivier} it follows that for all $M\in \N$ the map $\mathfrak{g}^{\ast,M}$ is an algebra homomorphism from $(\h^1_z,\ast)$ to $\C[[q]]$.

\begin{dfn}\label{def:gast}
For integers $s_1,\ldots,s_l\ge1$ and $M \in \N$, we define the $q$-series $\gstmrev_{s_1,\ldots,s_r}(q)\in\C[[q]]$ as the image of the word $w= z_{s_1}\dots z_{s_l} \in \h^1_z$ under the algebra homomorphism $(Z^{\ast}\otimes \mathfrak{g}^{\ast,M})\circ \Delta_H$:
\[ \gstmrev_{s_1,\ldots,s_l}(\tau) :=m\left( (Z^{\ast}\otimes \mathfrak{g}^{\ast,M})\circ \Delta_H \big( w \big) \right) \in \C[[q]]\,.\]
\end{dfn}
For $s_1,\dots,s_l \geq 2$ it is easy to see that the limit
\[ \gstrev_{s_1,\ldots,s_l}(\tau) :=  \lim_{M\to\infty}  \gstmrev_{s_1,\ldots,s_l}(\tau) \, \] 
exists and that we have 
\begin{prop}
For $s_1,\dots,s_l \geq 2$ we have $\grev_{s_1,\ldots,s_l} = \gstrev_{s_1,\ldots,s_l}=\gshrev_{s_1,\ldots,s_l}$.
\end{prop}
\begin{prf}
This follows since the construction above was exactly the one which appears in the calculation of the Fourier expansion of multiple Eisenstein series. See \cite{Bach} and \cite{BT} for details. 
\end{prf}

We now want to discuss whether for more general $s_1,\dots,s_l\in \N$ the limit of $\gstmrev_{s_1,\ldots,s_l}(\tau)$ as $M\to\infty$ exists. Since it is a finite sum of ordered sums of multitangent functions we can, by Theorem \ref{thm:olivier} iv), restrict to the case of ordered sums of monotangent functions and powers of $\pi$, i.e. we want to determine when the limit of 

\[ \sum_{0 < m_1 < \dots < m_l <M} f_1(m\tau) \dots f_l(m\tau) \] 
with $f_j(\tau) = \Psi_s(\tau)$ for some $s\in \N$ or $f_j(\tau) = 1$ exists. One easily checks that this exactly the case when $f_l(\tau)$ has no constant term, i.e. $f_l(\tau) \neq \Psi_1(\tau)$ and $f_l(\tau) \neq 1$. We deduce that therefore the limit of $\gstmrev_{s_1,\ldots,s_l}(\tau)$ as $M\to\infty$ exists when all $\Psi_{s_1,\dots,s_l}, \Psi_{s_2,\dots,s_l},\dots, \Psi_{s_l}$ have no constant term. Even though the Theorem \ref{thm:olivier} iv) just justifies this for the case all $s_j \geq 2$ we see, by using the explicit reductions to monotangents given in \cite{Boulliot}, that for low weights in fact the $\Psi_{1,\dots,1}(\tau)$ are the only multitangent functions with constant term. This question remains open but seems to be crucial in order to get a definition of $\gstrev$ for all admissible indices. The functions $g\tbinom{s_1,\dots,s_l}{r_1,\dots,r_l}$, i.e. the bi-brackets, will appear in $\gstrev_{s_1,\ldots,s_l}$ every time there is a $j<l$ with $s_j=1$ as we will see in the following examples:

\begin{ex} \begin{enumerate}[i)]
\item We are going to calculate $\gstrev_{2,1,2}$. For this we use the table at the end of \cite{Boulliot} where one can find that $\Psi_{2,1,2}(z) =\Psi_{1,2}(z) =\Psi_{2,1}(z)= 0$, therefore it is
\begin{align*}
 \gstmrev_{2,1,2}(\tau) &= \zetarev(2,1,2)^\ast +\zetarev(2,1)^\ast \sum_{0 < m_1 < M} \Psi_2(m_1 \tau)+ \zetarev(2)^\ast \sum_{0 < m_1 < m_2 < M} \Psi_1(m_1\tau) \Psi_2(m_2\tau)   \\
 &+ \sum_{0 < m_1 < m_2 < m_3 < M} \Psi_2(m_1 \tau) \Psi_1(m_2 \tau) \Psi_2(m_3 \tau) \,.
\end{align*}
Taking the limit $M \to \infty$ and using the explicit forms of $\Psi_k$ ($k \geq 1$) and $\zetarev(2,1)^\ast = -\zetarev(1,2)-\zetarev(3)=-2\zetarev(1,2)$ we obtain 
\begin{align*}
  \gstrev_{2,1,2}= &\lim_{M\to\infty} \gstmrev_{2,1,2} \\
  = &\zetarev(2,1,2)-2 \zetarev(1,2) g_2 + \zetarev(2)\left( g_{1,2} + \frac{1}{2} g\tbinom{2}{1} - \frac{(-2\pi i)}{2} g_{2} \right) \\
  &+ g_{2,1,2}+\frac{1}{2}\left( g\tbinom{2,2}{1,0}-g\tbinom{2,2}{0,1}-(-2\pi i)g_{2,2} \right)   \\
  = &\zetarev(2,1,2) -2 \zetarev(1,2) g^\sh_{2} + \zetarev(2) g^\sh_{1,2}  + g^\sh_{2,1,2}\,\\
  = &\gshrev_{2,1,2} \,.
 \end{align*}
Similarly one can prove that $\gshrev_{1,2} =  \gstrev_{1,2}$, $\gshrev_{1,2,2} = \gstrev_{1,2,2}$ and $\gshrev_{1,4} = \gstrev_{1,4}$ from which we obtain the following stuffle product in weight $5$:
\begin{equation}\label{eq:stufflew5}
 \gshrev_2 \cdot \gshrev_{1,2} = \gshrev_{2,1,2} + 2 \gshrev_{1,2,2} + \gshrev_{3,2} + \gshrev_{1,4}\,.
\end{equation}

\item There are $\gstrev_{s_1,\dots,s_l}$ that differ from $\gshrev_{s_1,\dots,s_l}$. For example it is
\begin{align*}
\gstrev_{1,1,2} &= \zetarev(1,1,2) -\frac{13}{2} \zetarev(2) g_2 - (-2\pi i) g_{1,2} + \frac{1}{2} g\tbinom{1,2}{0,1}-\frac{3}{8} (-2\pi i)g\tbinom{2}{1} +\frac{1}{4}g\tbinom{2}{2}  + g_{1,1,2}\,,\\
\gshrev_{1,1,2} &=  \zetarev(1,1,2) - 4\zetarev(2) g_2 - (-2\pi i) g_{1,2}+ \frac{1}{2} g\tbinom{1,2}{0,1} - \frac{3}{12} (-2\pi i)g\tbinom{2}{1}+\frac{1}{6} g\tbinom{2}{2} + g_{1,1,2} \,.\\
&\gshrev_{1,1,2}-\gstrev_{1,1,2}  = \frac{5}{2} g_2  +\frac{1}{8} (-2\pi i)g\tbinom{2}{1}-\frac{1}{12}g\tbinom{2}{2}  \neq 0 
\end{align*}
It is still an open question for which indices $s_1,\dots,s_l$ we have $\gshrev_{s_1,\dots,s_l} = \gstrev_{s_1,\dots,s_l}$. The author wants to address this question in upcoming projects.  
\end{enumerate}
\end{ex} 

We end this paper by a comparison of different version of the double shuffle relations. For this we write for words $u,v \in \h^1$, $ds(u,v) = u \shuffle v - u \ast v \in \h^1$, where the $\shuffle$ is again the shuffle product with respect to the alphabet $\{x,y\}$ and $\ast$ the stuffle product with respect to the alphabet $\{z_1,z_2,\dots\}$. Recall that $\h^0$ denotes the set of all admissible words and set $\h^2=\Q\langle \{z_2,z_3,\dots \} \rangle$ to be the span of all words in $\h^1$ with no $z_1$ occurring, i.e. the words for which the  multiple Eisenstein series $G$ exists.

With this we define the numbers $eds_k$ (extended double shuffle relations of weight $k$), $fds_k$ (finite double shuffle relations of weight $k$) and $rds_k$ (restricted finite double shuffle relations of weight $k$) by
\begin{align*}
eds_k &:= \dim_\Q \big\langle ds(u,v) \in \h^1 \mid |u|+|v|=k, \,\,  u \in \h^0 , v \in \h^0 \cup \{ y \}   \big\rangle_\Q \,,\\
fds_k &:= \dim_\Q \big\langle ds(u,v) \in \h^1 \mid |u|+|v|=k ,\,\, u,v \in \h^0   \big\rangle_\Q \,,\\
rds_k &:= \dim_\Q \big\langle ds(u,v) \in \h^1 \mid  |u|+|v|=k ,\,\, u,v \in \h^2   \big\rangle_\Q \,.
\end{align*}
For the number of admissible generators of weight $k$ which equals $2^{k-2}$ for $k>1$, i.e. words in $\h^0$, we write $gen_k$. By Theorem \ref{thm:restrds} we know that the number of relations between the $\gshrev$ of weight $k$ is at least $rds_k$. But these relations don't suffice to give all relations between (shuffle) regularized multiple Eisenstein series since some of the finite double shuffle relations which are not restricted are also fulfilled. The numbers $d_k$ and $d'_k$ defined by
\[ \sum_{k\geq 0} d_k X^k = \frac{1}{1-X^2-X^3} \,, \quad \sum_{k\geq 0} d'_k X^k = \frac{1-X^2+X^4}{1-2X^2-2X^3}  \,,\]  
 are the conjectured dimensions for $\MZ_k$ and $\grw_k(\MDA)$ (see \cite{BK} Remark 5.7) respectively. Since it is also conjectured that $eds_k$ is the number of all relations between MZV of weight $k$ one expects that $d_k = gen_k-eds_k$, which so far is not known.
  It was observed in \cite{BT} that up to weight $7$ the dimension of
    \[ \mathcal{E}_k = \big\langle \gshrev_{s_1,\ldots,s_l}(q) \mid k=s_1+\cdots+s_l,\ l\ge0,s_1,\dots,s_{l-1}\ge1,s_l\ge2 \big\rangle_{\Q}\]
 seems to be the same as the dimension of $\grw_k(\MDA)$, i.e. conjecturally $d'_k$. We therefore set $cds_k := gen_k - d'_k$ which gives the number of conjectured relations in $\grw_k(\MDA)$ which coincide with the number of relations in $\mathcal{E}_k$ up to weight $7$ due to the calculations of the authors in \cite{BT}. The following table gives an overview of these numbers up to weight $14$.

\begin{table}[H]\footnotesize
\captionsetup{width=0.8\textwidth}
\begin{center}
\begin{tabular} {|c|c|c|c|c|c|c|c|c|c|c|c|c|c|c|}\hline
$k$ 		   & 1	&	2 & 3 & 4	& 5 & 6 &  7 & 8 & 9 & 10 & 11 & 12 & 13 & 14  \\ \hline
$eds_k$ 	 & 0	&	0 & 1 & 3 & 6	& 14 & 29 &  60 & 123 & 249 & 503 & 1012 & 2032 & 4075  \\ \hline
$fds_k$ 	 & 0	&	0 & 0 & 1	& 2 & 7 &  16 & 40 & 92 & 200 & 429 & 902 & 1865 & 3832  \\ \hline
$cds_k:=gen_k-d'_k$ 	 & 0	&	0 & 0 & 1	& 2 & 6 &  14 & 32 & 72 & 156 & 336 & 712 & 1496 & 3120  \\ \hline
$rds_k$ 	 & 0	&	0 & 0 & 1	& 1 & 3 &  5 & 11 & 19 & 37 & 65 & 120 & 209 &  372 \\ \hline \hline
$gen_k$ 	 & 0	&	1 & 2 & 4	& 8 & 16 &  32 & 64 & 128 & 256 & 512 & 1024 & 2048 & 4096  \\ \hline
$d_k \overset{?}{=} gen_k-eds_k  $ 		  & 0	&	1 & 1 & 1	& 2 & 2 &  3 & 4 & 5 & 7 & 9 & 12 & 16 & 21  \\ \hline
$d'_k  $    & 0 & 1 & 2 & 3 & 6 & 10 & 18 & 32 & 56 & 100 & 176 & 312 &552 &976 	\\ \hline
$\dim \mathcal{E}_k \geq $    & 0 & 1 & 2 & 3 & 6 & 10 & 18 & ? & ? & ? & ? & ? &? &? 	\\ \hline
\end{tabular}
\caption{Comparison of the number of extended-, finite-, conjectured- and restricted-double shuffle relations.} 
\end{center}
\end{table} 
The last line give lower bounds of the dimension of the space $\mathcal{E}_k$ spanned by all admissible shuffle regularized multiple Eisenstein series of weight $k$ which are for $k\leq 5$ exact since we derived all relations up to this weight.

{\small
{\it E-mail:}\\\texttt{henrik.bachmann@uni-hamburg.de}\\

\noindent {\sc Fachbereich Mathematik (AZ)\\ Universit\"at Hamburg\\ Bundesstrasse 55\\ D-20146 Hamburg}}

\end{document}